%
%
%
%


\magnification=1200
\vsize=20.4cm \hsize=13.5cm
\parskip=4pt plus 2pt minus 1pt
\parindent=5.2mm

\font\twentybf = ec-lmbx10 at 20.736pt
-lmbx10 at 14.4pt
\font\twelvebf = ec-lmbx10 at 12pt

\font\twelvebsy=cmbsy10 at 12pt
\font\tenbsy=cmbsy10
\font\eightbsy=cmbsy8
\font\sevenbsy=cmbsy7
\font\sixbsy=cmbsy6
\font\fivebsy=cmbsy5

\font\tenCal=eusm10 at 10pt
\font\eightCal=eusm10 at 8pt
\font\sevenCal=eusm10 at 7pt
\font\sixCal=eusm10 at 6pt
\font\fiveCal=eusm10 at 5pt

\newfam\Calfam
  \textfont\Calfam=\tenCal
  \scriptfont\Calfam=\sevenCal
  \scriptscriptfont\Calfam=\fiveCal
\def\Cal{\fam\Calfam\tenCal}


\font\tenmsa = msam10
\font\sevenmsa = msam7
\font\fivemsa = msam5
\newfam\msafam
\textfont\msafam=\tenmsa
\scriptfont\msafam=\sevenmsa
\scriptscriptfont\msafam=\fivemsa

\font\tenmsb = msbm10
\font\sevenmsb = msbm7
\font\fivemsb = msbm5
\newfam\msbfam
\textfont\msbfam=\tenmsb
\scriptfont\msbfam=\sevenmsb
\scriptscriptfont\msbfam=\fivemsb

\font\teneuf = eufm10
\font\seveneuf = eufm7
\font\fiveeuf = eufm5
\newfam\euffam
\textfont\euffam=\teneuf
\scriptfont\euffam=\seveneuf
\scriptscriptfont\euffam=\fiveeuf

\catcode`\@=11

\def\Bbb#1{{\fam\msbfam#1}}
\def\hexnumber@#1{\ifcase#1 0\or1\or2\or3\or4\or5\or6\or7\or8\or9\or
 A\or B\or C\or D\or E\or F\fi}
\edef\msa@{\hexnumber@\msafam}
\edef\msb@{\hexnumber@\msbfam}
\mathchardef\square="0\msa@03
\mathchardef\smallsetminus="2\msb@72
\let\ssm=\smallsetminus
\mathchardef\compact="3\msa@62
\mathchardef\upharpoonright="3\msa@16
\let\restriction\upharpoonright
\mathchardef\subsetneq="3\msb@28
\mathchardef\supsetneq="3\msb@29
\catcode`\@=12


\def\ii{{\rm i}}
\def\bu{{\scriptstyle\bullet}}
\def\bul{{\scriptscriptstyle\bullet}}

\def\build#1^#2_#3{\mathrel{\mathop{\null#1}\limits^{#2}_{#3}}}
\def\buildo#1\over#2{\mathrel{\mathop{\null#2}\limits^{#1}}}
\def\buildu#1\under#2{\mathrel{\mathop{\null#2}\limits_{#1}}}

\font\twelvemib = lmmib10 at 12pt
\font\tenmib = lmmib10
\font\sevenmib = lmmib10 at 7pt

\font\eightrm = ec-lmr10 at 8pt
\font\eightbf = ec-lmbx10 at 8pt
\font\eightit = ec-lmri10 at 8pt
\font\eighttt = ec-lmtt10 at 8pt
\font\eighti = lmmi10 at 8pt
\font\eightsy = lmsy10 at 8pt
\font\eightmib = lmmib10 at 8pt
\font\eightex = lmex10 at 8pt

\font\eightmsa = msam8
\font\eightmsb = msbm8
\font\eighteuf = eufm8

\font\sixrm = ec-lmr10 at 6pt
\font\sixbf = ec-lmbx10 at 6pt
-lmri10 at 6pt
-lmtt10 at 6pt
\font\sixi = lmmi10 at 6pt
\font\sixsy = lmsy10 at 6pt
\font\sixmib = lmmib10 at 6pt

\font\sixmsa = msam6
\font\sixmsb = msbm6
\font\sixeuf = eufm6

\font\fiverm = ec-lmr10 at 5pt
\font\fivebf = ec-lmbx10 at 5pt
-lmri10 at 5pt
-lmtt10 at 5pt
\font\fivei = lmmi10 at 5pt
\font\fivesy = lmsy10 at 5pt

\font\fivemib = lmmib10 at 5pt

\font\fivemsa = msam5
\font\fivemsb = msbm5

\def\eightpoint{\def\rm{\fam0\eightrm}
   \textfont0=\eightrm \scriptfont0=\sixrm \scriptscriptfont0=\fiverm
   \textfont1=\eighti  \scriptfont1=\sixi  \scriptscriptfont1=\fivei
   \textfont2=\eightsy \scriptfont2=\sixsy \scriptscriptfont2=\fivesy
   \textfont3=\eightex \scriptfont3=\eightex\scriptscriptfont3=\eightex
   \def\it{\fam\itfam\eightit}%
   \textfont\itfam=\eightit
   \def\bf{\fam\bffam\eightbf}%
   \textfont\bffam=\eightbf \scriptfont\bffam=\sixbf
   \scriptscriptfont\bffam=\fivebf
   \def\tt{\fam\ttfam\eighttt}%
   \textfont\ttfam=\eighttt
   \textfont\msafam=\eightmsa \scriptfont\msafam=\sixmsa
   \textfont\msbfam=\eightmsb \scriptfont\msbfam=\sixmsb
   \textfont\euffam=\eighteuf \scriptfont\euffam=\sixeuf
   \textfont\Calfam=\eightCal
   \scriptfont\Calfam=\sixCal
   \normalbaselineskip=9.6pt
   \normalbaselines\rm}

\newskip\eightpointsurround
\eightpointsurround=6pt plus 3pt minus 2pt
\def\begineightpoint{\vskip\eightpointsurround \bgroup\eightpoint}
\def\endeightpoint{\vskip\eightpointsurround \egroup}

\def\twelvepointbf{%
 \textfont0=\twelvebf   \scriptfont0=\eightbf   \scriptscriptfont0=\sixbf
 \textfont1=\twelvemib  \scriptfont1=\eightmib \scriptscriptfont1=\sixmib
 \textfont2=\twelvebsy  \scriptfont2=\eightbsy  \scriptscriptfont2=\sixbsy
 \twelvebf
 \baselineskip=14.4pt}

\def\tenpointbf{%
 \textfont0=\tenbf   \scriptfont0=\sevenbf   \scriptscriptfont0=\fivebf
 \textfont1=\tenmib  \scriptfont1=\sevenmib  \scriptscriptfont1=\fivemib
 \textfont2=\tenbsy  \scriptfont2=\sevenbsy  \scriptscriptfont2=\fivebsy
 \tenbf
 \baselineskip=14.4pt}

\def\today{\ifcase\month\or
January\or February\or March\or April\or May\or June\or July\or August\or
September\or October\or November\or December\fi \space\number\day,
\number\year}

\catcode`\@=11
\newcount\@tempcnta \newcount\@tempcntb 
\def\timeofday{{%
\@tempcnta=\time \divide\@tempcnta by 60 \@tempcntb=\@tempcnta
\multiply\@tempcntb by -60 \advance\@tempcntb by \time
\ifnum\@tempcntb > 9 \number\@tempcnta:\number\@tempcntb
  \else\number\@tempcnta:0\number\@tempcntb\fi}}
\catcode`\@=12

\def\proof#1{\removelastskip\vskip4pt\noindent{\it #1.}}
\def\endproof{\vskip5pt plus 1pt minus 2pt}

\def\qed{\relax\ifmmode
            \hbox{$\square$}%
         \else
            {\unskip\nobreak\hfil
            \penalty50\hskip1em\null\nobreak\hfil\hbox{$\square$}%
            \parfillskip=0pt\finalhyphendemerits=0\endgraf}%
         \fi}
\def\bigsquare{{\kern-0.3ex\hbox{
\vrule height 1.7ex  width 0.093ex  depth 0ex\kern-0.093ex
\vrule height 1.8ex  width 1.7ex  depth -1.707ex\kern-0.093ex
\vrule height 1.7ex  width 0.093ex  depth 0ex\kern-1.65ex
\vrule height 0.093ex  width 1.6ex  depth 0ex}\kern0.3ex}}

\long\def\claim#1 #2\endclaim
{\removelastskip\vskip7pt plus 1pt minus 1pt\noindent{\bf#1.}
{\it\ignorespaces#2}\vskip\baselineskip}

\newdimen\plainitemindent\plainitemindent=18pt

\def\\{\hfil\break}

\def\bC{{\Bbb C}}
\def\bN{{\Bbb N}}

\def\bR{{\Bbb R}}

\def\cC{{\Cal C}}
\def\cI{{\Cal I}}
\def\cJ{{\Cal J}}
\def\cO{{\Cal O}}

\def\stimes{\mathop{\kern0.7pt
\vrule height 0.4pt depth 0pt width 5pt\kern-5pt
\vrule height 5.4pt depth -5pt width 5pt\kern-5pt
\vrule height 5.4pt depth 0pt width 0.4pt\kern4.6pt
\vrule height 5.4pt depth 0pt width 0.4pt\kern-6.5pt
\raise0.3pt\hbox{$\times$}\kern-0.7pt}\nolimits}

\catcode`@=11
\def\rightarrowfil{\m@th\mathord-\mkern-6mu%
  \cleaders\hbox{$\mkern-2mu\mathord-\mkern-2mu$}\hfill
  \mkern-6mu\mathord\rightarrow}
\catcode`@=12

\def\Ll{\langle\!\langle}
\def\Rr{\rangle\!\rangle}

\def\dbar{\overline\partial}
\def\ddbar{\partial\overline\partial}
\def\Ker{\mathop{\rm Ker}\nolimits}
\def\Id{\mathop{\rm Id}\nolimits}

\def\Re{\mathop{\rm Re}\nolimits}
\def\Im{\mathop{\rm Im}\nolimits}
\def\mod{\mathop{\rm mod}\nolimits}

\def\Supp{\mathop{\rm Supp}\nolimits}
\def\Jac{\mathop{\rm Jac}\nolimits}
\def\Hom{\mathop{\rm Hom}\nolimits}
\def\Dom{\mathop{\rm Dom}\nolimits}
\def\Nak{{\rm Nak}}

\def\reg{{\rm reg}}

\def\item#1{\par\parindent=\plainitemindent\noindent
\hangindent\parindent\hbox to\parindent{#1\hss}\ignorespaces}


\newcount\void \void=-1
\newcount\newtitle \newtitle=\void
\newbox\titlebox   \setbox\titlebox\hbox{\hfil}
\newbox\sectionbox \setbox\sectionbox\hbox{\hfil}
\def\folio{\ifnum\pageno=1 {\hfil}
           \else
           \ifodd\pageno
           {\eightpoint\hfil\copy\sectionbox\hfil\hfil\bf\number\pageno}\else
           {\eightpoint\bf\number\pageno\hfil\hfil\copy\titlebox\hfil}
           \fi\fi}
\footline={\hfil}
\headline={\folio}

\def\supereject{%
  \vfill\eject$\strut$\ifodd\pageno \else\newtitle=\void
  \headline={\hfil} \vfill\eject \headline={\folio} \fi}

\def\author#1{\noindent{\twelvebf #1}\vskip3pt}
\def\address#1{\noindent #1}

\def\titlerunning#1{\setbox\titlebox\hbox{\eightpoint\hskip-25mm\bf
  \abbreviatedauthorname, #1}}
\def\title#1{%
  \newtitle=\pageno
  \vskip5pt\noindent{\twentybf #1}\noindent
  \titlerunning{#1}}

\def\realbreak#1{\break\phantom{\kern #1}}
\def\fakebreak#1{ }
\def\realhfill{\hfill}
\def\fakehfill{ }
\def\realdash{-}
\def\fakedash{\kern-0.01pt}
\def\sectionrunning#1#2{\setbox\sectionbox\hbox{\eightpoint\bf #1. #2}}
\def\section#1#2{%
  \let\sectionbreak=\fakebreak
  \let\sectionhfill=\fakehfill
  \let\sectiondash=\fakedash
  \par\vskip0.75cm\penalty -100
  \let\sectionbreak=\realbreak
  \let\sectionhfill=\realhfill
  \let\sectiondash=\realdash
  \vbox{\baselineskip=14.4pt\noindent{{\twelvepointbf #1.\kern6pt #2}}}
  \sectionrunning{#1}{#2}
  \vskip5pt
  \penalty 500}

\def\nonumsectionrunning#1{\setbox\sectionbox\hbox{\eightpoint #1}}

\def\subsection#1#2{%
  \par\vskip0.2cm\penalty -100
  \vbox{\noindent{{\tenpointbf #1. #2}}}
  \vskip1pt
  \penalty 500}

\def\subsubsection#1#2{%
  \par\vskip0.1cm\penalty -100
  \noindent{{\tenpointbf #1.} {\it #2}}}

\def\cite#1{[#1]}

\def\Bibitem[#1]&#2&#3&#4&{\noindent\hangindent=\parindent\hangafter=1
{\bf [#1]} {\rm #2}, {\it #3}, {\rm #4.}\vskip2pt plus 0.5pt minus 0.5pt}

\def\setheadtitles#1{%
\setbox\sectionbox\hbox{\eightpoint #1}
\setbox\titlebox\hbox{\eightpoint #1}}

\def\twolines{\\\phantom{$~$}\hfill\\}


\def\abbreviatedauthorname{J.-P.\ Demailly}

\twolines\twolines
\title{Extension of holomorphic}
\title{functions defined on non}\vskip3pt
\title{reduced analytic subvarieties}\vskip15pt
\titlerunning{Extension of holomorphic functions}

\author{Jean-Pierre Demailly}
\address{Institut Fourier, Universit\'e de Grenoble-Alpes}
\footnote{*}{\eightrm This work is supported by the ERC grant ALKAGE}
\vskip2cm

\noindent
{\eightpoint
{\bf Abstract.} The goal of this contribution is to investigate $L^2$
extension properties for holomorphic sections of vector bundles satisfying
weak semi-positivity properties. Using techniques borrowed from recent
proofs of the Ohsawa-Takegoshi extension theorem, we obtain several
new surjectivity results for the restriction morphism to a non
necessarily reduced subvariety, provided the latter is defined as 
the zero variety of a multiplier ideal sheaf. These extension results
come with precise $L^2$ estimates and
(probably) optimal curvature conditions.\bigskip

\bigskip

\noindent
{\bf MSC Classification 2010.} 14C30, 14F05, 32C35
}
\vskip1.5cm

\section{1}{Introduction}

Bernhard Riemann is considered to be the founder of many branches of
geometry, especially by providing the foundations of modern
differential geometry [Rie54] and by introducing a deep geometric
approach in the theory of functions of one complex variable [Rie51,
Rie57]. The study of Riemannian manifolds has since become a central
theme of mathematics.  We are concerned here with the study of
holomorphic functions of several variables, and in this context,
Hermitian geometry is the relevant special case of Riemannian
manifolds -- as a matter of fact, Charles Hermite, born in 1822, was a
contemporary of Bernhard Riemann.

The more specific problem we are considering is the question whether a
holomorphic function $f$ defined on a subvariety $Y$ of a complex
manifold $X$ can be extended to the whole of~$X$. More generally,
given a holomorphic vector bundle $E$ on $X$, we are interested in the
existence of global extensions $F\in H^0(X,E)$ of a section $f\in
H^0(Y,E_{\restriction Y})$, assuming suitable convexity properties of
$X$ and $Y$, suitable growth conditions for $f$, and appropriate
curvature positivity hypotheses for the bundle $E$. This can also be
seen as an interpolation problem when $Y$ is not connected -- possibly
a discrete set.

In his general study of singular points and other questions of
function theory, Riemann obtained what is now known as the Riemann
extension theorem: a singularity of a holomorphic function at a point
is ``removable'' as soon as it is bounded -- or more generally if it
is $L^2$ in the sense of Lebesgue; the same is true for functions of
several variables defined on the complement of an analytic set. Using
the maximum principle, a consequence is that a holomorphic function
defined in the complement of an analytic set of codimension 2
automatically extends to the ambient manifold.

The extension-interpolation problem we are considering is strongly
related to Riemann's extension result because we allow here $Y$ to
have singularities, and $f$ must already have a local extension near
any singular point $x_0\in Y$; for this, some local growth condition
is needed in general, especially if $x_0$ is a non normal point. A
major advance in the general problem is the Ohsawa-Takegoshi $L^2$
extension theorem [OT87] (see also the subsequent series of papers
II--VI by T.~Ohsawa). It is remarkable that Bernhard Riemann already
anticipated in [Rie51] the use of $L^2$ estimates and the idea of
minimizing energy, even though his terminology was very different from
the one currently in use.

The goal of the present contribution is to generalize the Ohsawa-Takegoshi
$L^2$ extension theorem to the case where the sections are extended
from a non necessarily reduced subvariety, associated with an
arbitrary multiplier ideal sheaf. A similar idea had already been
considered in D.~Popovici's work [Pop05], but the present approach is
substantially more general. Since the extension theorems do not
require any strict positivity assumptions, we hope that they will be
useful to investigate further properties of linear systems and
pluricanonical systems on varieties that are not of general type. The
exposition is organized as follows: section~2 presents the main definitions
and results, section 3 (which is more standard and which the expert
reader way wish to skip) is devoted to recalling the required Kähler
identities and inequalities, section 4 elaborates on the concept of
jumping numbers for multiplier ideal sheaves, and section 5 explains
the technical details of the $L^2$ estimates and their proofs. We refer
to [Ber96], [BL14], [Blo13], [Che11], [DHP13], [GZ13,15], [Man93], [MV07], 
[Ohs88,94,95,01,03,05], [Pop05], [Var10] for related work on $L^2$ extension 
theorems.

The author adresses warm thanks to Mihai P\u{a}un and Xiangyu Zhou for
several stimu\-lating discussions around these questions.

\section{2}{Statement of the main extension results}

The Ohsawa-Takegoshi theorem addresses the following extension problem:
let $Y$ be a complex analytic submanifold of a complex manifold~$X\,$;
given a holomorphic function $f$ on $Y$ satisfying suitable $L^2$
conditions on $Y$, find a holomorphic extension $F$ of $f$ to $X$,
together with a good $L^2$ estimate for $F$ on~$X$. For this,
suitable pseudoconvexity and curvature conditions are needed.
We start with a few basic definitions in these directions. All
manifolds and complex spaces are assumed to be paracompact (and even
countable at infinity).

\claim{(2.1) Definition} A compact complex manifold $X$ will be said
to be {\rm weakly pseudoconvex} if it possesses a smooth
plurisubharmonic exhaustion function $\rho$ $($say $\rho:X\to\bR_+)$.
\smallskip\noindent
{\rm(Note: in [Nak73] and later work by Ohsawa, the terminology 
``weakly 1-complete'' is used instead of ``weakly pseudoconvex'')}.
\endclaim

\claim{(2.2) Definition} A function $\psi:X\to [-\infty,+\infty[$ on a
complex manifold $X$ is said to be {\rm quasi-plurisubharmonic}
$($quasi-psh$)$ if $\varphi$ is locally the sum of a psh function and
of a smooth function $($or equivalently, if $\ii\ddbar\psi$ is locally
bounded from below$)$. In addition, we say that $\psi$ has {\rm neat
analytic singularities} if every point $x_0\in X$ possesses an open
neighborhood $U$ on which $\psi$ can be written
$$\psi(z)=c\log\sum_{1\le j\le N}|g_k(z)|^2+w(z)$$
where $c\ge 0$, $g_k\in\cO_X(U)$ and $w\in\cC^\infty(U)$.
\endclaim

\claim{(2.3) Definition} If $\psi$ is a quasi-psh function on a
complex manifold $X$, the {\it multiplier ideal sheaf} $\cI(\psi)$ is
the coherent analytic subsheaf of $\cO_X$ defined by
$$
\cI(\psi)_x=\big\{f\in\cO_{X,x}\,;\;\exists U\ni x\,,\;
\int_U|f|^2e^{-\psi}d\lambda<+\infty\big\}
$$
where $U$ is an open coordinate neighborhood of $x$, and $d\lambda$
the standard Lebesgue measure in the corresponding open chart of $\bC^n$. We
say that the singularities of $\psi$ are {\rm log canonical} along the
zero variety $Y=V(\cI(\psi))$ if~
$\cI((1-\varepsilon)\psi)_{\restriction Y}=\cO_{X\restriction Y}$ for
every $\varepsilon>0$.
\endclaim

In case $\psi$ has log canonical singularities, it is easy to see that
$\cI(\psi)$ is a reduced ideal, i.e.\ that $Y=V(\cI(\psi))$ is a reduced
analytic subvariety of $X$. In fact if $f^p\in\cI(\psi)_x$ for some
integer $p\ge 2$, then $|f|^{2p}e^{-\psi}$ is locally integrable, hence, by
openness, we have also $|f|^{2p}e^{-(1+\varepsilon)\psi}$ for $\varepsilon>0$ 
small enough (see [GZ13] for a general result on openness that does not assume
anything on $\psi$ -- the present special case follows in fact directly from
Hironaka's theorm on resolution of singularities, in case $\psi$ has 
analytic singularities). Using the fact that 
$e^{-(1-\varepsilon)\psi}$ is integrable 
for every $\varepsilon\in{}]0,1[$, the Hölder inequality for the conjugate
exponents $1/p+1/q=1$ implies that 
$$|f|^2e^{-\psi}=\big(|f|^{2p}e^{-(1+\varepsilon)\psi}\big)^{1/p}
\big(e^{-\psi(1/q - \varepsilon/p)}\big)$$ is locally integrable,
hence $f\in\cI(\psi)_x$, as was to be shown.  If $\omega$ is a Kähler
metric on~$X$, we let $dV_{X,\omega}={1\over n!}\omega^n$ be the corresponding
K\"ahler volume element, $n=\dim X$. In case $\psi$ has log canonical
singularities along $Y=V(\cI(\psi))$, one can also associate in a
natural way a measure $dV_{Y^\circ,\omega}[\psi]$ on the set
$Y^\circ=Y_\reg$ of regular points of $Y$ as follows. If
$g\in\cC_c(Y^\circ)$ is a compactly supported continuous function on
$Y^\circ$ and $\widetilde g$ a compactly supported extension of $g$ to
$X$, we set
$$
\int_{Y^\circ}g\,dV_{Y^\circ,\omega}[\psi]=
\limsup_{t\to-\infty}\int_{\{x\in X\,,\;t<\psi(x)<t+1\}}
\widetilde ge^{-\psi}\,dV_{X,\omega}.\leqno(2.4)
$$
By Hironaka, it is easy to see that the limit does not depend on the
continuous extension~$\widetilde g$, and that one gets in this way a
measure with smooth positive density with respect to the Lebesgue
measure, at least on an (analytic) Zariski open set in $Y^\circ$ (cf.\
Proposition 4.5, in the special case $f=1$, $p=1$). 
In case $Y$ is a codimension $r$ subvariety of $X$
defined by an equation $\sigma(x)=0$ associated with a section $\sigma
\in H^0(X,S)$ of some hermitian vector bundle $(S,h_S)$ on $X$, and
assuming that $\sigma$ is generically transverse to zero along $Y$, it is
natural to take
$$\psi(z)=r\log|\sigma(z)|^2_{h_S}.\leqno(2.5)$$
One can then easily check that $dV_{Y^\circ,\omega}[\psi]$ is the
measure supported on $Y^\circ=Y_{\rm reg}$ such~that
$$
dV_{Y^\circ,\omega}[\psi]={2^{r+1}\pi^r\over (r-1)!}\,
{1\over|\Lambda^r(d\sigma)|^2_{\omega,h_S}}dV_{Y,\omega}\quad
\hbox{where}\quad
dV_{Y,\omega}= {1\over (n-r)!}\,\omega^{n-r}_{\restriction Y^\circ}.\leqno(2.6)
$$
For a quasi-psh function with log canonical singularities, 
$dV_{Y^\circ,\omega}[\psi]$ should thus be seen as some sort of
(inverse of) Jacobian determinant associated with the logarithmic 
singularities of~$\psi$. Finally, the following positive real function
will make an appearance in several of our final $L^2$ estimates~:
$$
\gamma(x) =\left\{\eqalign{
\rlap{\hbox{$\exp(-x/2)$}}\kern60pt &\hbox{if $x\ge 0$},\cr
\rlap{\hbox{$\displaystyle{1\over 1+x^2}$}}\kern60pt &\hbox{if $x\le 0$}.\cr}
\right.\leqno(2.7)
$$
The first generalized extension theorem we are interested in
is a variation of Theorem~4 in [Ohs01]. The first difference
is that we do not require any specific behavior of the quasi-psh function
defining the subvariety: any quasi-psh function with log canonical 
singularities will do; secondly, we do not want to make any assumption
that there exist negligible sets in the ambient manifold whose complements 
are Stein, because such an hypothesis need not be true
on a general compact K\"ahler manifold -- one of the targets of our study.
Recall that a hermitian tensor 
$$
\Theta = \ii\sum_{1\le j,k\le n,\,1\le\lambda,\mu\le r}
c_{jk\lambda\mu}dz_j\wedge d\overline z_k\otimes e_\lambda^*\otimes e_\mu
\in \cC^\infty(X,\Lambda^{1,1}T^*_X\otimes\Hom(E,E))
$$
is said to be {\it Nakano semi-positive} (resp.\ {\it positive}) if
the associated hermitian quadratic form
$$
H_\Theta(\tau)=\sum_{1\le j,k\le n,\,1\le\lambda,\mu\le r}c_{jk\lambda\mu}
\tau_{j,\lambda}\overline \tau_{k,\mu}
$$
is semi-positive (resp.\ positive) on non zero tensors
$\tau=\sum\tau_{j,\lambda}{\partial\over\partial z_j}\otimes e_\lambda\in
T_X\otimes E$. It~is said to be  {\it Griffiths semi-positive} (resp.\ 
{\it positive}) if $H_\Theta(\tau)\ge 0$ (resp.\ \hbox{$H_\Theta(\tau)>0$})
for all non zero decomposable tensors $\tau=\zeta\otimes v\in T_X\otimes E$.

\claim{(2.8) Theorem (general $L^2$ extension result for reduced subvarieties)}
Let $X$ be a weakly pseudoconvex Kähler manifold, and $\omega$ a
Kähler metric on~$X$.  Let $(E,h)$ be a holomorphic vector bundle 
equipped with a smooth hermitian metric $h$ on $X$, and let
$\psi:X\to[-\infty,+\infty[$ be a quasi-psh function on $X$ with neat
analytic singularities. Let $Y$ be the analytic subvariety of $X$ defined 
by $Y=V(\cI(\psi))$ and assume that $\psi$ has log canonical
singularities along $Y$, so that $Y$ is reduced. Finally, assume
that the Chern curvature tensor of $(E,h)$ is such that the sum
$$
\ii\Theta_{E,h}+\alpha\,\ii\ddbar\psi\otimes\Id_E
$$
is Nakano semipositive for
all $\alpha\in[1,1+\delta]$ and some $\delta>0$. Then for every
section $f\in H^0(Y^\circ,(K_X\otimes E)_{\restriction Y^\circ})$ on $Y^\circ=Y_\reg$
such that
$$
\int_{Y^\circ}|f|^2_{\omega,h}dV_{Y^\circ,\omega}[\psi]<+\infty,
$$
there exists an extension $F\in H^0(X,K_X\otimes E)$ whose restriction to
$Y^\circ$ is equal to $f$, such that
$$
\int_X\gamma(\delta\psi)\,|F|^2_{\omega,h}e^{-\psi}
dV_{X,\omega}\le{34\over\delta}\int_{Y^\circ}|f|^2_{\omega,h}dV_{Y^\circ,\omega}[\psi].
$$
\endclaim

\claim{(2.9) Remarks} {\rm (a) Although $|F|^2_{\omega,h}$ and
$dV_{X,\omega}$ both depend on $\omega$, it is easy to see that the
product $|F|^2_{\omega,h}dV_{X,\omega}$ actually does not depend on
$\omega$ when $F$ is a $(n,0)$-form. The same observation applies to
the product $|f|^2_{\omega,h}dV_{Y^\circ,\omega}[\psi]$, hence the
final $L^2$ estimate is in fact independent of $\omega$. Nevertheless,
the existence of a K\"ahler metric (and even of a complete K\"ahler
metric) is crucial in the proof, thanks to the techniques developped in
[AV65] and [Dem82].
\smallskip
\noindent
(b) By approximating non smooth plurisubharmonic 
weights with smooth ones, one can see that the above result still holds
when $E$ is a line bundle equipped with a singular hermitian metric
$h=e^{-\varphi}$. The curvature condition then reads
$$\ii\Theta_{E,h}+\alpha\,\ii\ddbar\psi=
\ii\,\ddbar(\varphi+t\psi)\ge 0,\qquad \alpha\in[1,1+\delta],$$
and should be understood in the sense of currents.\smallskip
\noindent
(c) The constant ${34\over\delta}$ given in the above $L^2$ inequality 
is not optimal. By exercising more care in the bounds, an optimal estimate
could probably be found by following the techniques of Blocki [Blo13] and
Guan-Zhou [GZ15], at the expense of replacing $\gamma(\delta\psi)$ with
a more complicated and less explicit function of $\delta$ and $\psi$.
Notice also that in the $L^2$ estimate, $\delta$ can be replaced by any
$\delta'\in{}]0,\delta]$.\qed}
\endclaim

We now turn ourselves to the case where non reduced multiplier ideal sheaves
and non reduced subvarieties are considered. This situation has already
been considered by D.~Popovici [Pop05] in the case of powers of a
reduced ideal, but we aim here at a much wider generality, which also
yields more natural assumptions. For $m\in\bR_+$, we consider the
multiplier ideal sheaf $\cI(m\psi)$ and the associated non necessarily reduced
subvariety $Y^{(m)}=V(\cI(m\psi))$, together with the structure sheaf
$\cO_{Y^{(m)}}=\cO_X/\cI(m\psi)$, the real number $m$ being viewed as some sort
of multiplicity -- the support $|Y^{(m)}|$ may increase with $m$, but certainly
stabilizes to set of poles $P=\psi^{-1}(-\infty)$ for $m$ large enough. 
We assume the existence of a discrete sequence of positive numbers
$$0=m_0<m_1<m_2<\ldots<m_p<\ldots$$
such that $\cI(m\psi)=\cI(m_p\psi)$ for $m\in [m_p,m_{p+1}[$ (with of course
$\cI(m_0\psi)=\cO_X)$; they
are called the {\it jumping numbers} of~$\psi$. The existence of a discrete
sequence of jumping numbers is automatic if $X$ is compact. In general, it
still holds on every relatively compact open subset 
$$X_c:=\{x\in X\,,\;\rho(x)<c\}\compact X,$$ 
but requires some some of uniform behaviour of singularities at infinity
in the noncompact case. We are interested in extending a holomorphic section 
$$
f\in H^0(Y^{(m_p)},\cO_{Y^{(m_p)}}(K_X\otimes E_{\restriction Y^{(m_p)}})
:=H^0(Y^{(m_p)},\cO_X(K_X\otimes_\bC E)\otimes_{\cO_X}\cO_X/\cI(m_p\psi)).
$$
[Later on, we usually omit to specify the rings over which tensor products
are taken, as they are implicit from the nature of objects under consideration].
The results are easier to state in case one takes a nilpotent section of
the form
$$
f\in H^0(Y^{(m_p)},\cO_X(K_X\otimes E)\otimes\cI(m_{p-1}\psi)/\cI(m_p\psi)).
$$
Then $\cI(m_{p-1}\psi)/\cI(m_p\psi))$ is actually a coherent sheaf, and its
support is a reduced subvariety $Z_p$ of $Y^{(m_p)}$ (see Lemma~4.2). Therefore 
$\cI(m_{p-1}\psi)/\cI(m_p\psi))$ can be seen as a vector bundle
over a Zariski open set $Z_p^\circ\subset Z_p$. We can mimic formula (2.4)
and define some sort of infinitesimal ``$m_p$-jet'' $L^2$ norm
$|J^{m_p}f|^2_{\omega,h}\,dV_{Z_p^\circ,\omega}[\psi]$ (a purely formal
notation), as the measure on $Z_p^\circ$ defined by
$$
\int_{Z_p^\circ}g\,|J^{m_p}f|^2_{\omega,h}\,dV_{Z_p^\circ,\omega}[\psi]=
\limsup_{t\to-\infty}\int_{\{x\in X\,,\;t<\psi(x)<t+1\}}
\widetilde g\,|\widetilde f|^2_{\omega,h}e^{-m_p\psi}\,dV_{X,\omega}\leqno(2.10)
$$
for any $g\in\cC_c(Z_p^\circ)$, where $\widetilde g\in\cC_c(X)$ is a 
continuous extension of $g$ and $\widetilde f$ a smooth
extension of $f$ on $X$ such that $\smash{\widetilde f} - f\in
\cI(m_p\psi)\otimes_{\cO_X}\cC^\infty$ (this measure again has a
smooth positive density on a Zariski open set in $Z_p^\circ$, and does
not depend on the choices of $\smash{\widetilde f}$ and $\widetilde
g$, see Prop.~4.5). We extend the measure as being $0$ on 
$\smash{Y^{(m_p)}_{\rm red}}\ssm Z_p$, since $\cI(m_{p-1}\psi)/\cI(m_p\psi))$ 
has support
in~$Z_p^\circ\subset Z_p$. In this context, we introduce the following
natural definition.

\claim{(2.11) Definition} We define the {\rm restricted multiplied ideal sheaf}
$$
\cI'(m_{p-1}\psi)\subset\cI(m_{p-1}\psi)
$$
to be the set of germs $F\in\cI(m_{p-1}\psi)_x\subset\cO_{X,x}$ such that 
there exists a neighborhood $U$ of $x$ satisfying
$$
\int_{Y^{(m_p)}\cap U}|J^{m_p}F|^2_{\omega,h}\,
dV_{Y^{(m_p)},\omega}[\psi]<+\infty.
$$
This is a coherent ideal sheaf that contains $\cI(m_p\psi)$.
Both of the inclusions
$$
\cI(m_p\psi)\subset \cI'(m_{p-1}\psi)\subset\cI(m_{p-1}\psi)
$$
can be strict $($even for $p=1)$.
\endclaim

The proof is given in Section~4 (Proposition 4.5).  One of the
geometric consequences is the following ``quantitative'' surjectivity
statement, which is the analogue of Theorem~2.8 for the case when the
first non trivial jumping number $m_1$ is replaced by a higher jumping
number~$m_p$.

\claim{(2.12) Theorem} With the above notation and in the general 
setting of Theorem~$2.8$ $($but without the hypothesis 
that the quasi-psh function $\psi$ has log canonical singularities$)$, let
$0=m_0<m_1<m_2<\ldots<m_p<\ldots$ be the jumping numbers of~$\psi$. Assume that
$$
\ii\Theta_{E,h}+\alpha\,\ii\ddbar\psi\otimes\Id_E\ge_\Nak 0
$$
is Nakano semipositive for all $\alpha\in[m_p,m_p+\delta]$, for some
$\delta>0$.\smallskip
\item{\rm(a)} Let
$$
f\in H^0(Y^{(m_p)},\cO_X(K_X\otimes E)\otimes\cI'(m_{p-1}\psi)/\cI(m_p\psi))
$$
be a section such that
$$
\int_{Y^{(m_p)}}|J^{m_p}f|^2_{\omega,h}\,dV_{Y^{(m_p)},\omega}[\psi]<+\infty.
$$
Then there exists a global section
$$
F\in H^0(X,\cO_X(K_X\otimes E)\otimes\cI'(m_{p-1}\psi))
$$
which maps to $f$ under the morphism $\cI'(m_{p-1}\psi)\to
\cI(m_{p-1}\psi)/\cI(m_p\psi)$, such that
$$
\int_{X}\gamma(\delta\psi)\,|F|^2_{\omega,h}\,e^{-m_p\psi}dV_{X,\omega}[\psi]
\le {34\over\delta}\int_{Y^{(m_p)}}|J^{m_p}f|^2_{\omega,h}\,dV_{Y^{(m_p)},\omega}[\psi].
$$
\item{\rm(b)} The restriction morphism
$$
H^0(X,\cO_X(K_X\otimes E)\otimes\cI'(m_{p-1}\psi))\to 
H^0(Y^{(m_p)},\cO_X(K_X\otimes E)\otimes\cI'(m_{p-1}\psi)/\cI(m_p\psi))
$$
is surjective.\vskip0pt
\endclaim

If $E$ is a line bundle and $h$ a singular hermitian metric on $E$, a
similar result can be obtained by approximating $h$. However, the
$L^2$ estimates require to take into account the multiplier ideal
sheaf of~$h$, and we get the following result.

\claim{(2.13) Theorem} Let $(X,\omega)$ be a weakly pseudoconvex 
K\"ahler manifold,
$\psi$ a quasi-psh function with neat analytic singularities and $E$ a
holomorphic line bundle equipped with a singular hermitian metric $h$.
Let $0=m_0<m_1<m_2<\ldots<m_p<\ldots$ be the jumping numbers for the family of 
multiplier ideal sheaves $\cI(he^{-m\psi})$, $m\in\bR^+$. Assume that
$$
\ii\Theta_{E,h}+\alpha\,\ii\ddbar\psi\ge 0\quad\hbox{in the sense of currents},
$$
for all $\alpha\in[m_p,m_p+\delta]$ and $\delta>0$ small enough.
Define $\cI'(he^{-m_\ell\psi})$ by the same $L^2$ convergence
property as in the case when $h$ is non singular. Then the restriction morphism
$$
\eqalign{
H^0(X,\cO_X(K_X\otimes E)&\otimes\cI'(he^{-m_{p-1}\psi}))\cr
&\to H^0(Y^{(m_p)},\cO_X(K_X\otimes E)\otimes\cI'(he^{-m_{p-1}\psi})/
\cI(he^{-m_p\psi}))\cr}
$$
is surjective. Moreover, the $L^2$ estimate {\rm (2.12)~(b)} still holds
in this situation.
\endclaim

One of the strengths of the above theorems lies in the fact that no
strict curvature hypothesis on $\ii\Theta_{E,h}$ is ever needed. On
the other hand, if we assume a strictly positive lower bound
$$
\ii\Theta_{E,h}+m_p\,\ii\ddbar\psi\ge \varepsilon\omega>0,
$$
then the Nadel vanishing theorem [Nad89] implies \hbox{$H^1(X,
\cO_X(K_X{\otimes}E)\otimes\cI(he^{-m_p\psi}))=0$}, and we immediately
get stronger surjectivity statements by considering the relevant
cohomology exact sequence, e.g.\ that
$$
H^0(X,\cO_X(K_X\otimes E)\otimes\cI(he^{-m_\ell\psi}))
\to H^0(Y^{(m_p)},\cO_X(K_X\otimes E)\otimes\cI(he^{-m_\ell\psi})/
\cI(he^{-m_p\psi}))
$$
is surjective for all $\ell<p$. In case $X$ is compact, it turns out that 
this qualitative
surjectivity property still holds true with a semi-positivity
assumption only. In view of [DHP13], this probably has interesting
applications to algebraic geometry which we intend to discuss in
a future work:

\claim{(2.14) Theorem} Let $(X,\omega)$ be a {\rm compact}
K\"ahler manifold,
$\psi$ a quasi-psh function with neat analytic singularities and $E$ a
holomorphic vector bundle equipped with a hermitian metric $h$.
Let $0=m_0<m_1<m_2<\ldots<m_p<\ldots$ be the jumping numbers for the family of 
multiplier ideal sheaves $\cI(he^{-m\psi})$, $m\in\bR^+$. 
\smallskip
\item{\rm(a)} The rank of $E$ being arbitrary, assume that $h$ is 
smooth and that for some index
$\ell=0,1,\ldots,p-1$ and all $\alpha\in[m_{\ell+1},m_p+\delta]$, we have
$$
\ii\Theta_{E,h}+\alpha\,\ii\ddbar\psi\otimes\Id_E\ge 0\quad\hbox{in the sense of Nakano},
$$
for $\delta>0$ small enough. Then the restriction morphism
$$
H^0(X,\cO_X(K_X\otimes E)\otimes\cI(e^{-m_\ell\psi}))
\to H^0(Y^{(m_p)},\cO_X(K_X\otimes E)\otimes\cI(e^{-m_\ell\psi})/\cI(e^{-m_p\psi}))
$$
is surjective.
\smallskip
\item{\rm(b)} Assume that $E$ is a line bundle equipped with a singular
metric $h$ and that for some
$\ell=0,1,\ldots,p-1$ and all $\alpha\in[m_{\ell+1},m_p+\delta]$, we have
$$
\ii\Theta_{E,h}+\alpha\,\ii\ddbar\psi\ge 0\quad\hbox{in the sense of currents},
$$
for $\delta>0$ small enough. Then the restriction morphism
$$
\eqalign{
H^0(X,\cO_X(K_X\otimes E)&\otimes\cI(he^{-m_\ell\psi}))\cr
&\to H^0(Y^{(m_p)},\cO_X(K_X\otimes E)\otimes\cI(he^{-m_\ell\psi})/
\cI(he^{-m_p\psi}))\cr}
$$
is surjective.\vskip0pt
\endclaim

\claim{(2.15) Question} {\rm It would be interesting to know whether
Theorem~2.14 can be strengthened by suitable quantitative
$L^2$ estimates.  The main difficulty is already to define the norm
of jets when there is more than one jump number involved. Some sort of
``Cauchy inequality'' for jets would be needed in order to derive the
successive jet norms from a known global $L^2$ estimate for a holomorphic
section defined on the whole of~$X$. We do not know how to proceed
further at this point.}
\endclaim

\section{3}{Fundamental estimates of K\"ahler geometry}

\subsection{3.A}{Basic set-up}

We refer to [Gri66,69] and [Dem-X] for general background results on
the geometry of K\"ahler manifolds and hermitian bundles. Let $X$ be a 
complex $n$-dimensional manifold equipped with a smooth Hermitian metric
$$
\omega=\sum_{1\le j,k\le n}\omega_{jk}(z)\,dz_j\otimes d\overline z_k.
$$
As usual, it is convenient to view $\omega$ rather as a real
$(1,1)$-form $\omega=i\sum\omega_{jk}(z)\, dz_j\wedge d\overline
z_k$. The metric $\omega$ is said to be {\it K\"ahler} if $d\omega=0$,
i.e.\ $\omega$ also defines a symplectic structure.  It can be easily
shown that $\omega$ is K\"ahler if and only if there are holomorphic
coordinates $(z_1,\,\ldots\,,z_n)$ centered at any point $x_0\in X$
such that the matrix of coefficients $(\omega_{jk})$ is tangent to
identity at order $2$, i.e.\
$$\omega_{jk}(z)=\delta_{jk}+O(|z|^2)\quad\hbox{at $x_0$}.\leqno(3.1)$$
Now, let $(E,h)$ is a Hermitian vector bundle over $X$. Given a smooth
$(p,q)$-form  $u$ on $X$ with values in $E$, that is, a section of
$\cC^\infty(X,\Lambda^{p,q}T^*_X\otimes E)$, we consider the global $L^2$ norm
$$\|u\|^2 = \int_M |u(x)|^2dV_{X,\omega}(x)\leqno(3.2)$$
where $|u(x)|=|u(x)|_{\omega,h}$ is the pointwise Hermitian norm of $u(x)$ 
in the tensor product $\Lambda^{p,q}T^*_X\otimes E$ and 
$dV_{X,\omega}={\omega^n\over n!}$ the Hermitian volume form on~$X\,$; 
for simplicity, we will usually 
omit  the dependence on the metrics in the notation of $|u(x)|_{\omega,h}$.
We denote by $\Ll\bul,\bul\Rr$ the inner product of the Hilbert space
$L^2(X,\Lambda^{p,q}T^*_X\otimes E)$ of $L^2$ sections for the norm
$\|\bul\|$. Let $D=D_{E,h}$ be the Chern connection of $(E,h)$, that is,
the unique connection $D=D^{1,0}+D^{0,1}=D'+D''$ that makes $h$ parallel, 
for which the $(0,1)$ component $D''$ coincides with $\dbar$.
It follows in particular that $D^{\prime 2}=0$, $D^{\prime\prime 2}=0$ and
$D^2=D'D''+D''D'$. The {\it Chern curvature tensor} of $(E,h)$ is 
the $(1,1)$-form
$$\Theta_{E,h}\in \cC^\infty(X,\Lambda^{1,1}T^*_X\otimes\Hom(E,E))\leqno(3.3)$$
such that $D^2u=(D'D''+D''D')u=\Theta_{E,h}\wedge u$.
Next, the complex Laplace-Beltrami operators are defined by
$$\Delta'=D'D^{\prime*}+D^{\prime*}D',\qquad
\Delta''=D''D^{\prime\prime*}+D^{\prime\prime*}D''\leqno(3.4)$$
where $P^*$ denotes the formal adjoint of a differential operator 
$$P:\cC^\infty(X,\Lambda^{p,q}T^*_X\otimes E)\to
\cC^\infty(X,\Lambda^{p+r,q+s}T^*_X\otimes E),\qquad
\deg P=r+s$$
with respect to the corresponding inner products $\Ll\bul,\bul\Rr$. If 
$u$ is a compactly supported section, we get in particular
$$
\Ll \Delta' u,u\Rr=\|D'u\|^2+\|D^{\prime *}u\|^2\ge 0,\qquad
\Ll \Delta'' u,u\Rr=\|D''u\|^2+\|D^{\prime\prime *}u\|^2\ge 0.
$$

\subsection{3.B}{Bochner-Kodaira-Nakano identity}
Great simplifications occur in the operators of Hermitian geometry
when the metric $\omega$ is K\"ahler. In fact, if we use normal 
coordinates at a point $x_0$ (cf.\ $(3.1)$), and a local holomorphic frame
$(e_\lambda)_{1\le\lambda\le r}$ of $E$ such that $De_\lambda(x_0)=0$, it
is not difficult to see that all order 1 operators $D'$, $D''$ and their
adjoints $D^{\prime*}$, $D^{\prime\prime*}$ admit at $x_0$
the same expansion as the analogous operators obtained when all Hermitian
metrics on $X$ or $E$ are constant. From this, the basic commutation
relations of K\"ahler geometry can be checked. If $A,B$ are differential
operators acting on the algebra $\cC^\infty(X,\Lambda^{\bu,\bu}T^*_X\otimes E)$,
their graded commutator (or graded Lie bracket) is $[A,B]= AB - (-1)^{ab} BA$
where $a,b$ are the degrees of $A$ and $B$ respectively. If $C$ is another
endomorphism of degree $c$, the following purely formal {\it Jacobi identity}
holds:
$$(-1)^{ca} \big[A,[B,C]\big] + (-1)^{ab}\big[B,[C,A]\big]
+ (-1)^{bc} \big[C,[A,B]\big] = 0.$$

\claim{(3.5) Fundamental K\"ahler identities} 
Let $(X,\omega)$ be a K\"ahler
manifold and let $L$ be the Lefschetz operator defined by $Lu=\omega\wedge u$ 
and $\Lambda=L^*$ its adjoint, acting on $E$-valued forms. The following
identities hold for the Chern connection $D=D'+D''$ on $E$ and the associated
complex Laplace operators $\Delta'$ and $\Delta''$.
\smallskip
\item{\rm(a)} Basic commutation relations 
$$[D^{\prime\prime*},L]=\ii D',~~[\Lambda,D'']=-\ii D^{\prime*},~~
[D^{\prime*},L]=-\ii D'',~~[\Lambda,D']=\ii D^{\prime\prime*}.$$
\smallskip
\item{\rm(b)} Bochner-Kodaira-Nakano identity 
{\rm([Boc48], [Kod53a,b], [AN54], [Nak55])}
$$\Delta''=\Delta'+[\ii\Theta_{E,h},\Lambda].$$
\endclaim

\proof{Idea of proof} (a) The first step is to check the identity
$[d^{\prime\prime*},L]=\ii d'$ for constant metrics on $X=\bC^n$ and the
trivial bundle $E=X\times\bC$, by a brute force calculation. All three 
other identities follow by taking conjugates or adjoints. The case of 
variable metrics follows by looking at Taylor expansions up to order~1.

\noindent
(b) The last equality in (a) yields $D^{\prime\prime*}
=-\ii[\Lambda,D']$, hence
$$\Delta''=[D'',D^{\prime\prime*}]=-\ii[D'',\big[\Lambda,D']\big].$$
By the Jacobi identity we get
$$\big[D'',[\Lambda,D']\big]=\big[\Lambda,[D',D'']]+\big[D',[D'',\Lambda]\big]
=[\Lambda,\Theta_{E,h}]+\ii[D',D^{\prime*}],$$
taking into account that $[D',D'']=D^2=\Theta_{E,h}$. The formula follows.\qed
\endproof

One important (well known) fact is that the curvature term
$[\ii\Theta_{E,h},\Lambda]$ operates as a hermitian (semi-)positive operator
on $L^2(X,\Lambda^{n,q}T^*_X\otimes E$ as soon as $\ii\,\Theta_{E,h}$ is
Nakano (semi-)positive.
\medskip

\subsection{3.C}{The twisted a priori inequality of Ohsawa and Takegoshi}

\noindent The main a priori inequality that we are going to use is a 
simplified (and slightly extended) version of the original
Ohsawa-Takegoshi a priori inequality [OT87, Ohs88], along the lines
proposed by Manivel [Man93] and Ohsawa [Ohs01]. Such inequalities were
originally introduced in the work of Donnelly-Fefferman [DF83] and
Donnelly-Xavier [DX84]. The main idea is to introduce a
modified Bochner-Kodaira-Nakano inequality. Although it has become
classical in this context, we reproduce here briefly the calculations 
for completeness, and also for the sake of fixing the notation.

\claim{(3.6) Lemma {\rm (Ohsawa [Ohs01])}} Let $E$ be a Hermitian
vector bundle on a complex manifold $X$ equipped with a K\"ahler
metric~$\omega$. Let $\eta,\,\lambda>0$ be smooth functions
on~$X$. Then for every form $u\in\cC^\infty_c(X,\Lambda^{p,q} T^*_X\otimes E)$
with compact support we have
$$
\eqalign{
\|(\eta+\lambda)^{{1\over2}}D^{\prime\prime*}u\|^2
&{}+\|\eta^{{1\over2}}D''u\|^2+\|\lambda^{{1\over2}}D'u\|^2+
2\|\lambda^{-{1\over2}}d'\eta\wedge u\|^2\cr
 &\ge \Ll[\eta\,\ii\Theta_E-\ii\,d'd''\eta-\ii\lambda^{-1}d'\eta\wedge
d''\eta,\Lambda]u,u\Rr.\cr}
$$
\endclaim
\proof{Proof} We consider the ``twisted'' Laplace-Beltrami operators
$$
\eqalign{
D'\eta D^{\prime*} +D^{\prime*}\eta  D'
& =  \eta [D',D^{\prime*}]+[D',\eta]D^{\prime*}+
[D^{\prime*},\eta]D' \cr
& = \eta\Delta'+(d'\eta)D^{\prime*}-(d'\eta)^*D',\cr
D''\eta D^{\prime\prime*} +D^{\prime\prime*}\eta  D''
& = \eta [D'',D^{\prime\prime*}]+[D'',\eta]D^{\prime\prime*}+
[D^{\prime\prime*},\eta]D'' \cr
& = \eta\Delta''+(d''\eta)D^{\prime\prime*}-(d''\eta)^*D'',\cr}
$$
where $\eta$, $(d'\eta)$, $(d''\eta)$ are abbreviated notations for the
multiplication operators\break
$\eta\bu$, $(d'\eta)\wedge\bu$, $(d''\eta)\wedge\bu$.
By subtracting the above equalities and taking into account the
Bochner-Kodaira-Nakano identity $\Delta''-\Delta'=[\ii\Theta_E,\Lambda]$,
we get
$$
\leqalignno{
D''&{}\eta D^{\prime\prime*}+D^{\prime\prime*}\eta  D'' -
D'\eta D^{\prime*}  -D^{\prime*}\eta  D' \cr
&{}=\eta [\ii\Theta_E,\Lambda]+(d''\eta)D^{\prime\prime*}
-(d''\eta)^* D''+(d'\eta)^* D'-(d'\eta)D^{\prime*}.
&(3.7)\cr}
$$
Moreover, the Jacobi identity yields
$$
[D'',[d'\eta,\Lambda ]]-[d'\eta,[\Lambda,D'']]+[\Lambda,[D'',d'\eta]]=0,
$$
whilst $[\Lambda,D'']=-\ii D^{\prime*}$ by the basic commutation
relations~3.5~(a). A straightforward computation shows that
$[D'',d'\eta]=-(d'd''\eta)$ and $[d'\eta,\Lambda ]=\ii(d''\eta)^*$.
Therefore we get
$$
\ii[D'',(d''\eta)^*]+\ii[d'\eta,D^{\prime*}]-[\Lambda,(d'd''\eta)]=0,
$$
that is,
$$
[\ii\,d'd''\eta,\Lambda]=
[D'',(d''\eta)^*]+[D^{\prime*},d'\eta]=
D''(d''\eta)^*+(d''\eta)^* D''+
D^{\prime*}(d'\eta)+(d'\eta)D^{\prime*}.
$$
After adding this to (3.7), we find
$$
\eqalign{
D''\eta D^{\prime\prime*}&{}+D^{\prime\prime*}\eta  D'' -
D'\eta D^{\prime*}  -D^{\prime*}\eta  D'
+[\ii\,d'd''\eta,\Lambda]
\cr
&{}=\eta [\ii\Theta_E,\Lambda]+(d''\eta)D^{\prime\prime*}+
D''(d''\eta)^*+(d'\eta)^* D'+D^{\prime*}(d'\eta).
\cr}
$$
We apply this identity to a form $u\in\cC_c^\infty(X,\Lambda^{p,q}T^*_X\otimes E)$
and take the inner bracket with~$u$. Then
$$
\Ll (D''\eta D^{\prime\prime*})u,u\Rr=\Ll\eta D^{\prime\prime*}u,
D^{\prime\prime*}u\Rr=\|\eta^{{1\over2}}D^{\prime\prime*}u\|^2,
$$
and likewise for the other similar terms. The above equalities
imply
$$
\eqalign{
&\|\eta^{{1\over 2}}D^{\prime\prime*}u\|^2+\|\eta^{{1\over 2}}D''u\|^2
-\|\eta^{{1\over2}}D'u\|^2-\|\eta^{{1\over2}}D^{\prime*}u\|^2\cr
&\qquad=\Ll[\eta\,\ii\Theta_E-\ii\,d'd''\eta,\Lambda]u,u\Rr
+2\Re\,\Ll D^{\prime\prime*}u,(d''\eta)^* u\Rr
+2\Re\,\Ll D'u,d'\eta\wedge u\Rr.\cr}
$$
By neglecting the negative terms $-\|\eta^{{1\over2}}D'u\|^2
-\|\eta^{{1\over2}}D^{\prime*}u\|^2$ and adding the squares
$$
\eqalign{
\|\lambda^{{1\over 2}}D^{\prime\prime*}u\|^2+
2\Re\,\Ll D^{\prime\prime*}u,(d''\eta)^* u\Rr+
\|\lambda^{-{1\over 2}}(d''\eta)^* u\|^2
&\ge 0,\cr
\|\lambda^{{1\over 2}}D'u\|^2+
2\Re\,\Ll D'u,d'\eta\wedge u\Rr+
\|\lambda^{-{1\over 2}}d'\eta\wedge u\|^2
&\ge0\cr}
$$
we get
$$
\leqalignno{
\|\eta^{{1\over 2}}
D^{\prime\prime*}u\|^2&{}+\|\lambda^{{1\over 2}}D''u\|^2
+\|\lambda^{{1\over 2}}D'u\|^2+\|\lambda^{-{1\over 2}}d'\eta\wedge u\|^2
+\|\lambda^{-{1\over 2}}(d''\eta)^* u\|^2\cr
&\ge\Ll[\eta\,\ii\Theta_E-\ii\,d'd''\eta,\Lambda]u,u\Rr.&(3.8)\cr}
$$
Finally, we use the identity $a^*=i[\overline a,\Lambda]$ for any $(1,0)$-form
$a$ to get
$$
(d'\eta)^*(d'\eta)-(d''\eta)(d''\eta)^*
=\ii[d''\eta,\Lambda](d'\eta)+\ii(d''\eta)[d'\eta,\Lambda]
=[\ii d''\eta\wedge d'\eta,\Lambda],
$$
which implies
$$
\|\lambda^{-{1\over 2}}d'\eta\wedge u\|^2-\|\lambda^{-{1\over2}}(d''\eta)^* u\|^2=-\Ll[\ii\lambda^{-1}d'\eta\wedge d''\eta,\Lambda]u,u\Rr.\leqno(3.9)
$$
The inequality asserted in Lemma~3.6 follows by adding (3.8) and (3.9).\qed
\endproof

In the special case of $(n,q)$-forms, the forms $D'u$ and $d'\eta\wedge u$
are of bidegree $(n+1,q)$, hence the estimate takes the simpler form
$$
\|(\eta+\lambda)^{{1\over2}}D^{\prime\prime*}u\|^2{}
+\|\eta^{{1\over2}}D''u\|^2\ge
\Ll[\eta\,\ii\Theta_E-\ii\,d'd''\eta
-\ii\lambda^{-1}\,d'\eta\wedge d''\eta,\Lambda]u,u\Rr.
\leqno(3.10)
$$

\subsection{3.D}{Abstract $L^2$ existence theorem for solutions 
of $\dbar$-equations}

\noindent Using standard arguments from functional analysis -- actually just basic pro\-per\-ties of Hilbert spaces -- the a~priori inequality (3.10) implies
a powerful $L^2$ existence theorem for solutions of $\dbar$-equations.

\claim{(3.11) Proposition} Let $X$ be a complete K\"ahler manifold
equipped with a $($non ne\-ces\-sarily complete$)$ K\"ahler metric
$\omega$, and let $(E,h)$ be a Hermitian vector bundle over~$X$.
Assume that there are smooth and bounded functions $\eta,\,\lambda>0$
on $X$ such that the $($Hermitian$)$ curvature operator
$$
B=B^{n,q}_{E,h,\omega,\eta,\lambda}=
[\eta\,\ii\Theta_{E,h}-\ii\,d'd''\eta-\ii\lambda^{-1}d'\eta\wedge d''\eta,
\Lambda_\omega]
$$
is positive definite everywhere on $\Lambda^{n,q}T^*_X\otimes E$, for
some $q\ge 1$. Then for every form $g\in L^2(X,\Lambda^{n,q}T^*_X\otimes E)$
such that $D''g=0$ and $\int_X\langle B^{-1}g,g\rangle
\,dV_{X,\omega}<+\infty$, there exists $f\in
L^2(X,\Lambda^{n,q-1}T^*_X\otimes E)$ such that  $D''f=g$ and
$$
\int_X(\eta+\lambda)^{-1}
|f|^2\,dV_{X,\omega}\le \int_X\langle B^{-1}g,g\rangle\,dV_{X,\omega}.
$$
\endclaim
\proof{Proof} Assume first that $\omega$ is complete K\"ahler metric.
Let $v\in L^2(X,\Lambda^{n,q}T^*_X\otimes E)$, and
$v=v_1+v_2\in(\Ker D'')\oplus(\Ker D'')^\perp$ the decomposition of $v$
with respect to the closed subspace $\Ker D''$ and its orthogonal.
Since $g\in\Ker D''$, The Cauchy-Schwarz inequality yields
$$
|\Ll g,v\Rr|^2=|\Ll g,v_1\Rr|^2=|\Ll B^{-{1\over 2}}g,B^{{1\over 2}}v_1\Rr|^2\le
\int_X\langle B^{-1}g,g\rangle \,dV_{X,\omega}
\int_X\langle Bv_1,v_1\rangle \,dV_{X,\omega},
$$
and provided that $v\in\Dom D^{\prime\prime*}$, we find $v_2\in
(\Ker D'')^\perp\subset(\Im D'')^\perp=\Ker D^{\prime\prime *}$, and so
$D''v_1=0$, $D^{\prime\prime*}v_2=0$, whence
$$
\int_X\langle Bv_1,v_1\rangle \,dV_{X,\omega}\le
\|(\eta+\lambda)^{{1\over 2}}D^{\prime\prime*}v_1\|^2+
\|\eta^{{1\over 2}}D''v_1\|^2
=\|(\eta+\lambda)^{{1\over 2}}D^{\prime\prime*}v\|^2.
$$
Combining both inequalities, we obtain
$$
|\Ll g,v\Rr|^2
\le \Big(\int_X\langle B^{-1}g,g\rangle \,dV_{X,\omega}\Big)
\|(\eta+\lambda)^{{1\over 2}}D^{\prime\prime*}v\|^2.
$$
The Hahn-Banch theorem applied to the linear form
$(\eta+\lambda)^{{1\over 2}}D^{\prime\prime*}v\mapsto \Ll v,g\Rr$
implies the existence of an element $w\in L^2(X,\Lambda^{n,q}T^*_X
\otimes E)$ such that
$$
\eqalign{
\|w\|^2&\le\int_X\langle B^{-1}g,g\rangle \,dV_{X,\omega}\qquad
\hbox{and}\cr
\Ll v,g\Rr&=\Ll(\eta+\lambda)^{{1\over2}}
D^{\prime\prime*}v,w\Rr\qquad
\forall g\in\Dom D''\cap\Dom D^{\prime\prime*}.\cr}
$$
It follows that $f=(\eta+\lambda)^{{1\over2}}w$ satisfies
$D''f=g$ as well as the desired $L^2$ estimate. If $\omega$ is not
complete, we set $\omega_\varepsilon=\omega+\varepsilon\widehat\omega$ with
some complete K\"ahler metric~$\widehat\omega$. The final conclusion is
then obtained by passing to the limit and using a monotonicity argument
(the integrals are easily shown to be monotonic with respect 
to~$\varepsilon$, see [Dem82]).\qed
\endproof

We need also a variant of the $L^2$-estimate, so as to obtain
approximate solutions with weaker requirements on the data$\,$:

\claim{(3.12) Proposition} With the notation of $3.11$, assume that 
$B+\varepsilon I>0$ for some $\varepsilon>0$ $($so that $B$ can be just
semi-positive or even slightly negative; here $I$ is the identity 
endomorphism$)$. Given a section
$g\in L^2(X,\Lambda^{n,q}T^*_X\otimes E)$ such that $D''g=0$ and
$$M(\varepsilon):=
\int_X\langle (B+\varepsilon I)^{-1}g,g\rangle\,dV_{X,\omega}<+\infty,$$
there exists an approximate solution
\hbox{$f_\varepsilon\in L^2(X,\Lambda^{n,q-1}T^*_X \otimes E)$} and a correcting
term $g_\varepsilon\in L^2(X,\Lambda^{n,q}T^*_X \otimes E)$ such that
$D''f_\varepsilon=g-g_\varepsilon$ and
$$
\int_X(\eta+\lambda)^{-1}|f_\varepsilon|^2\,dV_{X,\omega}+
{1\over \varepsilon}\int_X|g_\varepsilon|^2\,dV_{X,\omega}\le M(\varepsilon).
$$
If $g$ is smooth, then $f_\varepsilon$ and $g_\varepsilon$ can be taken smooth.
\endclaim

\proof{Proof} The arguments are almost unchanged, we rely instead on
the estimates
$$
|\Ll g,v_1\Rr|^2\le
\int_X\langle(B+\varepsilon I)^{-1}g,g\rangle \,dV_{X,\omega}
\int_X\langle(B+\varepsilon I)v_1,v_1\rangle \,dV_{X,\omega},
$$
and
$$
\int_X\langle(B+\varepsilon I)v_1,v_1\rangle \,dV_{X,\omega}\le
\|(\eta+\lambda)^{{1\over 2}}D^{\prime\prime*}v\|^2+
\varepsilon\|v\|^2.
$$
This gives a pair $(w_\varepsilon,w'_\varepsilon)$ such that 
$\Vert w_\varepsilon\Vert^2+\Vert w'_\varepsilon\Vert^2\le 
M(\varepsilon)$ and
$$
\Ll v,g\Rr = \Ll (\eta+\lambda)^{{1\over 2}}D^{\prime\prime*}v,w_\varepsilon\Rr+
\Ll\varepsilon^{1/2}v,w'_\varepsilon\Rr\quad
\hbox{for all $v\in\Dom D^{\prime\prime *}$},
\leqno(3.13)
$$
hence $f_\varepsilon=(\eta+\lambda)^{{1\over 2}}w_\varepsilon$ is the expected 
approximate solution with error term $g_\varepsilon=\varepsilon^{1/2}w'_\varepsilon$.
By (3.13), we do get $D''f_\varepsilon+g_\varepsilon=g$, and the 
expected $L^2$ estimates hold as well. In fact one can take 
$$f_\varepsilon=(\eta+\lambda)
D^{\prime\prime *}\square_\varepsilon^{-1}g\quad\hbox{and}\quad
g_\varepsilon=\varepsilon\,\square_\varepsilon^{-1}g
$$
where $\square_\varepsilon=D''(\eta+\lambda)D^{\prime\prime *}+D^{\prime\prime *}
(\eta+\lambda)D''+\varepsilon I$ is an invertible self-adjoint 
elliptic operator. Then $f_\varepsilon$ and $g_\varepsilon$ are smooth.\qed

\section{4}{Openness of multiplier ideal sheaves
and jumping numbers}

Let $X$ be complex manifold and $\varphi,\,\psi$ quasi-psh functions
on $X$. To every $m\in\bR_+$ we associate the multiplier ideal sheaf
$\cI(\varphi+m\psi)\subset\cO_X$. By Nadel [Nad89], this is a coherent
ideal sheaf, and $\cI(\varphi+m'\psi)\subset \cI(\varphi+m\psi)$ for
$m'\ge m$. The recent result of Guan-Zhou [GZ13] implies that one
has in fact $\cI(\varphi+(m+\alpha)\psi)=\cI(\varphi+m\psi)$ for every
$m\in\bR_+$ and $\alpha\in[0,\alpha_0(m)[$ sufficiently small. From this we
conclude without any restriction that there exists a discrete
sequence of numbers
$$0=m_0<m_1<m_2<\ldots<m_p<\ldots\leqno(4.1)$$
such that $\cI(\varphi+m\psi)=\cI(\varphi+m_p\psi)$ for $m\in[m_p,m_{p+1}[$ and
$$
\cI(\varphi)=\cI(\varphi+m_0\psi)\supsetneq\cI(\varphi+m_1\psi)\supsetneq
\ldots\supsetneq\cI(\varphi+m_l\psi)\supsetneq\ldots~.
$$
If $\psi$ is smooth, we have of course $m_1=+\infty$ already, and the sequence
stops there; in the sequel we assume that $\psi$ has non empty logarithmic 
poles to avoid this trivial situation -- the sequence $m_p$ is then
infinite since the multiplicities of germs of functions in 
$\cI(\varphi+m\psi)_x$ tend to infinity at every point $x$
where the Lelong number $\nu(\psi,x)$ is positive.

\claim{(4.2) Lemma} For every $p>0$, the ideal $\cJ_p\subsetneq\cO_X$
of germs of holomorphic functions $h$ such that
$h\,\cI(\varphi+m_{p-1}\psi)\subset\cI(\varphi+m_p\psi)$ is reduced,
i.e.\ $\sqrt{\cJ_p}=\cJ_p$. Moreover $\cJ_p$ contains
$\sqrt{\cI(m_p\psi)}$.
\endclaim

\proof{Proof} Let $h^k\in\cJ_{p,x}$ for some exponent $k\ge 2$. 
Pick $f\in\cI(\varphi+m_{p-1}\psi)_x$. By definition of the jumping 
numbers, $|f|^2e^{-\varphi-(m_p-\varepsilon)\psi}$
is integrable near $x$ for every $\varepsilon>0$. Since 
$h^kf\in\cI(\varphi+m_p\psi)$, the openness property shows that
$|h^kf|^2e^{-\varphi-(m_p+\delta)\psi}$ is integrable for some $\delta>0$.
For a suitable neighborhood $U$ of $x$ and $\varepsilon>0$ smaller than
$\delta/k$, the H\"older inequality applied with the measure
$d\mu=|f|^2e^{-\varphi-m_p\psi}d\lambda$ and the functions
$v=|h|^2e^{-\varepsilon\psi}$, $w=e^{\varepsilon\psi}$ for the conjugate
exponents $1/k+1/\ell=1$ implies
$$
\eqalign{
\int_U|hf|^2&e^{-\varphi-m_p\psi}d\lambda
=\int_Uvwd\mu\le\cr
&\Big(\int_U|h|^{2k}|f|^2e^{-\varphi-(m_p+k\varepsilon)\psi}d\lambda\Big)^{1/k}
\Big(\int_U|f|^2e^{-\varphi-(m_p-\ell\varepsilon)\psi}d\lambda
\Big)^{1/\ell}<+\infty,\cr}
$$
hence $hf\in\cI(\varphi+m_p\psi)$.
Since this is true for every $f\in\cI(\varphi+m_{p-1}\psi)$ we conclude that
$h\in\cJ_{p,x}$, thus $\sqrt{\cJ_p}=\cJ_p$. The last assertion is equivalent
to $\cJ_p\supset\cI(m_p\psi)$ and follows similarly from the inequality
$$
\eqalign{
\int_U|hf|^2&e^{-\varphi-m_p\psi}d\lambda\le\cr
&\Big(\int_U|h|^{2k}e^{-(m_p+k\varepsilon/\ell)\psi}d\lambda\Big)^{1/k}
\Big(\int_U|f|^{2\ell}e^{-\ell\varphi-(m_p-\varepsilon)\psi}d\lambda
\Big)^{1/\ell}<+\infty\,;\cr}
$$
we fix here $\varepsilon>0$ so small that 
$\cI(m_p\psi)=\cI((m_p+\varepsilon)\psi)$ and, by openness, $\ell$ sufficiently
close to $1$ to make the last integral convergent whenever
$f\in\cI(\varphi+m_{p-1}\psi)=\cI(\varphi+(m_p-\varepsilon)\psi)$.\qed

A consequence of Lemma 4.2 is that the zero variety $Z_p=V(\cJ_p)$ is a reduced
subvariety of $Y^{(m_p)}=V(\cI(m_p\psi))$ and that the quotient sheaf
$\cI(\varphi+m_{p-1}\psi)/\cI(\varphi+m_p\psi)$ is a coherent sheaf over
$\cO_{Z_p}=\cO_X/\cJ_p$. Therefore $\cI(\varphi+m_{p-1}\psi)/\cI(\varphi+m_p\psi)$
can be seen as a vector bundle on some Zariski open set
$Z_p^\circ\subset Z_p\subset \smash{Y^{(m_p)}_{\rm red}}$.

In the sequel, a case of special interest is when $\psi$ has analytic 
singularities, that is, every point $x_0\in X$ possesses an open
neighborhood $V\subset X$ on which $\psi$ can be written
$$\psi(z)=c\log\sum_{1\le j\le N}|g_k(z)|^2+u(z)$$
where $c\ge 0$, $g_k\in\cO_X(V)$ and $u\in\cC^\infty(V)$.  The integrability
of $|f|^2e^{-m\psi}$ on a coordinate neighborhood $V$ then means that
$$
\int_V{|f(z)|^2\over |g(z)|^{2m}}\,d\lambda(z)<+\infty.\leqno(4.3)
$$
By Hironaka [Hir64], there exists a principalization of the ideal 
$\cJ=(g_k)\subset\cO_X$, that is, a modification $\mu:\smash{\widehat X}
\to X$ such that $\mu^*\cJ=(g_k\circ\mu)=\cO_{\widehat X}(-\Delta)$ where
$\Delta=\sum a_k\Delta_k$ is a simple normal crossing divisor on $\smash{\widehat X}$.
We can also assume that the Jacobian $\Jac(\mu)$ has a zero divisor
$B=\sum b_k\Delta_k$ contained in the exceptional divisor. After a change
of variables $z=\mu(w)$ and a use of local coordinates where $\Delta_k=\{w_k=0\}$
and $\mu^*\cJ$ is the principal ideal generated by the monomial 
$w^a=\prod w_k^{a_k}$, we  see that (4.3) is equivalent to the convergence of
$$
\int_{\mu^{-1}(V)}{|f(\mu(w))|^2|\Jac(\mu)|^2\over 
|g(\mu(w)|^{2mc}}\,d\lambda(w),
$$
which can be expressed locally as the convergence of
$$
\int_{\mu^{-1}(V)}{|f(\mu(w))|^2|w^b|^2\over 
|w^a|^{2mc}}\,d\lambda(w).
$$
For this, the condition is that $f\circ\mu(w)$ be divisible by $w^s$ with
$s_k=\lfloor mca_k-b_k\rfloor_+$. In other words, the multipler ideal sheaves
$\cI(m\psi)$ are given by the direct image formula
$$
\cI(m\psi)=\mu_*\cO_{\widehat X}\big(-\sum_k\,
\lfloor mca_k-b_k\rfloor_+\Delta_k\big).\leqno(4.4)
$$
The jumps can only occur when $m$ is equal to one
of the values ${b_k+N\over ca_k}$, $N\in\bN$, which form a discrete subset 
of~$\bR_+$. 

\claim{(4.5) Proposition} Let $0=m_0<m_1<\ldots<m_p$ be the jumping
numbers of the quasi-psh function $\psi$, which is assumed to have
neat analytic singularities. Let $\ell$ a local holomorphic generator of
$K_X\otimes E$ at a point~$x_0\in X$, $E$ being equipped with a smooth 
hermitian metric $h$, let $Z_p=\Supp(\cI(m_{p-1}\psi)/\cI(m_p\psi))$, 
and take a germ $f\in\cI(m_{p-1}\psi)_{x_0}$.
\smallskip
\item{\rm(a)} The measure
$|J^{m_p}(f\ell)|^2_{\omega,h}\,dV_{Z_p^\circ,\omega}[\psi]$ defined by $(2.10)$ 
has a smooth
positive density with respect to the Lebesgue measure on a Zariski
open set $Z_p^\circ$ of~$Z_p$.
\smallskip
\item{\rm(b)} The sheaf $\cI'(m_{p-1}\psi)$ of germs 
$F\in\cI(m_{p-1}\psi)_x\subset\cO_{X,x}$ such that 
there exists a neighborhood $U$ of $x$ satisfying
$$
\int_{Y^{(m_p)}\cap U}|J^{m_p}(F\ell)|^2_{\omega,h}\,
dV_{Y^{(m_p)},\omega}[\psi]<+\infty
$$
is a coherent ideal sheaf such that
$$
\cI(m_p\psi)\subset \cI'(m_{p-1}\psi)\subset\cI(m_{p-1}\psi).
$$
Both of the inclusions can be strict $($even for $p=1)$.
\smallskip
\item{\rm(c)} The function $(1+|\psi|)^{-(n+1)}\,
|f\ell|^2_{\omega,h}e^{-m_p\psi}$ is locally integrable at~$x_0$.
\endclaim

\proof{Proof} (a) As above, we use a principalization of the
singularities of $\psi$ and apply formula~(4.4).  Over a generic point
$x_0\in Z_p$, the component of $Z_p$ containing $x_0$ is dominated by
exactly one of the divisors $\Delta_k$, and the jump number $m_p$ is
such that $m_pca_k-b_k=N$ for some integer~$N$.  The previous jump
number $m_{p-1}$ is then given by $m_{p-1}ca_k-b_k=N-1$.  On a
suitable coordinate chart of the blow-up $\mu:\smash{\widehat X}\to X$, let us
write $\psi\circ\mu(w)=c\log|w^a|^2+u(w)$ where $u$ is smooth, and let
$w^b=0$ be the zero divisor of $\Jac(\mu)$.  By definition, the
measure $|J^{m_p}f|^2_{\omega,h}\,dV_{Z_p^\circ,\omega}[\psi]$ is the
direct image of measures defined upstairs as the limit
$$
g\in\cC_c(Z_p,\bR)\mapsto \limsup_{t\to-\infty}\int_{t<c\log|w^a|^2+u(w)<t+1}
e^{-m_pu}\,\mu^*\widetilde g\;\beta(w)\,
|\widetilde v(w)|^2{e^{-\varphi\circ\mu}\over |w_k|^2}\,d\lambda(w).
$$
Here $\varphi$ is the weight of the metric $h$ on $E$, 
$\widetilde v(w)$ is the holomorphic function representing the
section $\mu^*(f\ell)/(w^{m_pca-b}/w_k)$ (the denominator divisor 
cancels with the numerator by construction), and $\beta(w)$ is
a smooth positive weight arising from the change of variable formula,
given by $\mu^*dV_{X,\omega}/|w^b|^2$ [one would still have to take
into account a partition of unity on the various coordinate charts
covering the fibers of $\mu$, but we will avoid this technicality for
the simplicity of notation]. Let us denote
$w=(w',w_k)\in\bC^{n-1}\times\bC$ and $d\lambda(w)=d\lambda(w')\,\lambda(w_k)$
the Lebesgue measure on $\bC^n$. At a generic point $w$ where
$w_j\ne 0$ for $j\ne k$, the domain of integration is of the form
$t(w')<ca_k\log|w_k|^2<t'(w')+1$. It is easy to check that the limsup measure
is a limit, equal to
$$
g\in\cC_c(Z_p,\bR)\mapsto \pi ca_k\,\big(e^{-m_pu}\,\mu^*\widetilde g\; 
|\widetilde v(w)|^2\beta(w)e^{-\varphi\circ\mu}
\big)_{\restriction w_k=0}\,d\lambda(w').
$$
We then have to integrate the right hand side measure over the fibers
of $\mu$ to get the density of this measure along $Z_p$. Since
$\mu$ can be taken to be a composition of blow-ups with smooth centers, 
the fibration has (upstairs) a locally trivial product structure over 
a Zariski open set $Z_p^\circ\subset Z_p$. Smooth local vector fields on
$Z_p^\circ$ can be lifted to smooth vector fields in the corresponding
chart of $\smash{\widehat X}$, and we conclude by differentiating
under the integral sign that the density downstairs on $Z_p$ 
is generically smooth.\smallskip

\noindent
(b) One typical example is given by 
$\psi(w)=\log|w_1|^2|w_2|^2$. Let us put $\Delta_k=\hbox{$\{w_k=0\}$}$, $k=1,2$. Then
the jumping numbers are $m_p=p\in\bN$, and we get
$$
\cI(m_p\psi)=\cO_X(-p(\Delta_1+\Delta_2)).
$$
For any $r_0>0$ fixed, we have
$$
\eqalign{
\int_{|w_1|<r_0,\,|w_2|<r_0,\,e^t<|w_1|^2|w_2|^2<e^{t+1}}
{1\over |w_1|^2|w_2|^2}\,&d\lambda(w_1)\,d\lambda(w_2)\cr
&\ge
\pi\int_{e^{t+1}/r_0^2<|w_2|^2<r_0^2}
{1\over |w_2|^2}\,d\lambda(w_2)\cr}
$$
and the limit as $t\to-\infty$ is easily seen to be infinite. Therefore
$\cI'(m_0\psi)=\cI_{\Delta_1\cap\Delta_2}$, and we do have
$$
\cO_X(-(\Delta_1+\Delta_2))=
\cI(m_0\psi)\subsetneq \cI'(m_0\psi)\subsetneq\cI(m_0\psi)=\cO_X
$$
in this case. In general, with our notation, it is easy to see that
$\cI'(m_{p-1}\psi)$ is given as the direct image
$$
\cI'(m_{p-1}\psi)=\mu_*\Big(\cO_{\widehat X}\big(-\sum_k\,
\lfloor m_{p-1}ca_k-b_k\rfloor_+\Delta_k\big)\otimes\cI_R\Big)
\leqno(4.6)
$$
where $R=\bigcup\Delta_\ell\cap\Delta_{\ell'}$ is the union of pairwise 
intersections of divisors $\Delta_\ell$ for which
$m_pca_\ell-b_\ell=m_pca_k-b_k\,({}=m_{p-1}ca_k-b_k+1)$, $k$ being
one of the indices achieving the values of $m_p$, $m_{p-1}$ at the given
point $x\in Z_p$. This is a coherent ideal sheaf by the direct image 
theorem.\smallskip

\noindent
(c) Near a point where $\psi\circ\mu(w)$ has singularities along normal
crossing divisors $w_j=0$, $1\le j\le k$, it is sufficient to show 
that integrals of the form
$$
I_k=\int_{|w_1|<1/2,\ldots,|w_k|<1/2}
{\big(-\log|w_1|^2\ldots|w_k|^2\big)^{-(k+1)}
\over |w_1|^2\ldots|w_k|^2}\,d\lambda(w_1)\ldots d\lambda(w_k)
$$
are convergent. This is easily done by induction on~$k$, by using
a change of variable $w_k'=w_1\ldots w_k$, and by applying the
Fubini formula together with an integration in $w_k'$.
We~then get $I_k\le {4\pi\over k}I_{k-1}$ (notice that 
$|w_k'|<|w_1|\ldots|w_{k-1}|$), and $I_1=\pi/(2\log 2)$, thus 
all $I_k$ are finite.\qed

\section{5}{Proof of the $L^2$ extension theorems}

Unless otherwise specified, $X$ denotes a weakly pseudoconvex complex
$n$-dimensional manifold equipped with a (non necessarily complete)
 K\"ahler metric~$\omega$, $\rho:X\to[0,+\infty[$ a smooth psh exhaustion 
on~$X$, $(E,h)$ a smooth hermitian holomorphic vector bundle, and 
$\psi:X\to[-\infty,+\infty[$ a quasi-psh function with neat analytic
singularities. Before given the technical details of the proofs, we 
start with a rather simple observation.

\claim{(5.1) Observation} {\rm Let $\mu:\widehat X\to X$ be a proper
modification. Assume that $\smash{\widehat X}$ is equipped with a K\"ahler
metric $\widehat\omega$.
\smallskip
\item{\rm(a)} For every $m\ge 0$, there is an isomorphism
$$
\mu^*:H^0(X,\cO(K_X\otimes E)\otimes\cI(m\psi))\to
H^0(X,\cO(K_{\widehat X}\otimes \mu^*E)\otimes\cI(m\psi\circ\mu))
$$
whose inverse is the direct image morphism $\mu_*$.
\smallskip
\item{\rm(b)} For any holomorphic section $F$ of $E$, the $L^2$ norm 
$$\int_{\widehat X}|F\circ \mu|^2e^{-m\psi\circ\mu}dV_{
\widehat X,\widehat\omega}$$
coincides with $\int_X|F|^2e^{-m\psi}dV_{X,\omega}$.
\smallskip
\item{\rm(c)} On the regular part of the subvariety $Y^{(m_p)}=
V(\cI(m_p\psi))$, for any
$$f\in H^0(Y,\cO_X(K_X\otimes E)\otimes\cI(m_{p-1}\psi)/\cI(m_p\psi)),$$
the area measure $|f|^2_{\omega,h}dV_{Y^{(m_p)},\omega}$ $($which is independent 
of $\omega)$ is the direct image of its
counterpart defined by $\widehat\psi=\psi\circ\mu$ and 
$\widehat f=\mu^*f$ on the strict transform of $Y^{(m_p)}$, i.e.\ the
union of components of
$\smash{\widehat Y}^{(m_p)}=\mu^{-1}(Y^{(m_p)})$ that have a
dominant projection to a component of $Y^{(m_p)}$.\qed\vskip0pt}
\endclaim

The proof of (c) is immediate by a change of variable $z=\mu(w)$ in the
integrals $\int_{\{t<\psi<t+1\}}...$ and by passing to the limits. It follows 
from the observation and the discussion of section~4 that after blowing up
the proof of our theorems can be reduced to the case where $\psi$ has 
divisorial singularities along a normal crossing divisor.\medskip

\noindent{\bf Proof of Theorem~2.8.} This will be only a mild
generalization of the techniques used in [Ohs01], with a technical
complication due to the fact that we do not assume any
Steinness of complements of negligible sets. With the notation of 
Theorem~2.8, let $f\in\smash{
H^0(Y^\circ,(K_X\otimes E)_{\restriction Y^\circ})}$.  We view $f$ as an
$E$-valued $(n,0)$-form defined over~$Y^\circ$ and apply Proposition~3.11
after replacing the metric $h$ of $E$ by the singular metric
$h_\psi=h^{-\psi}$, whose curvature is
$$
\ii\,\Theta_{E,h_\psi}=\ii\,\Theta_{E,h}+\ii\,d'd''\psi.
$$
In order to avoid the singularities, we shrink $X$ to a relatively compact
weakly pseudoconvex domain $X_c=\{\rho<c\}$, and work on $X_c\ssm Y$ instead
of $X$. In fact, we know by [Dem82, Theorem~1.5] that $X_c\ssm Y$ is 
complete K\"ahler. Let us first assume that $Y$ is non singular, i.e.\
$Y^\circ=Y$, and that
$\psi\le 0\,$; of course, when $X$ is compact, one can always
subtract a constant to $\psi$ to achieve $\psi\le 0$, but there could exist
non compact situations when it is interesting to take $\psi$ unbounded
from above. In any case, we claim that there exists a smooth section
$$
\widetilde f\in \cC^\infty(X,\Lambda^{n,0}T^*_X\otimes E)
$$
such that
{\parindent=0cm
\smallskip
(5.2) $\widetilde f$ coincides with $f$ on $Y$,
\smallskip
(5.3) $D''\widetilde f=0$ at every point of~$Y$,
\smallskip
(5.4) $|D''\widetilde f|_{\omega,h}^2e^{-\psi}$ is locally integrable near $Y$.
\smallskip}

\noindent
For this, consider a locally finite covering of $Y$ by coordinates 
patches $U_j\subset X$ biholomorphic to polydiscs, on which
$E_{\restriction U_j}$ is trivial, and such that either 
$Y\cap U_j=\emptyset$ or
$$
Y\cap U_j=\{z\in U_j\,;\;z_1=\cdots=z_r=0\}.
$$
We can find holomorphic sections
$$\widetilde f_j\in\cC^\infty(X,\Lambda^{n,q}T^*_X\otimes E)$$
which extend $f_{\restriction Y\cap U_j}$ [such sections can be obtained simply 
by viewing functions of the form $g(z_{r+1},\ldots,z_n)$ as independent of
$z_1,\ldots,z_r$ on each polydisc]. For some partition of unity $(\chi_j)$
subordinate to $(U_j)$, we then set
$\smash{\widetilde f:=\sum_j \chi_j\widetilde f_j}$. Clearly (5.2)
$\smash{\widetilde f}_{\restriction Y}=f$ holds. Since we have
$$
D''\widetilde f=D''(\widetilde f-\widetilde f_k)=
D''\Big(\sum_j\chi_j(\widetilde f_j-\widetilde f_k)\Big)=
\sum_jd''\chi_j\wedge(\widetilde f_j-\widetilde f_k),
$$
properties (5.3) and (5.4) also hold, as $\widetilde f_j-\widetilde f_k=0$
on $Y$ and $\cI(\psi)=\cI_Y$ is reduced
by our assumption that $\psi$ has log canonical singularities.
The main idea is to apply Proposition~3.11 to solve the equation
$$
D''u_t=v_t:=D''(\theta(\psi-t)\,\widetilde f),\qquad t\in{}]-\infty,-1],
\leqno(5.5)
$$
where $\theta:[-\infty,+\infty[{}\to [0,1]$ is a smooth non increasing function such that $\theta(\tau)=1$ for $\tau\in{}]-\infty,\varepsilon/3]$, 
$\theta(\tau)=0$ for $\tau\in[\varepsilon/3,+\infty[$ and $|\theta'|\le 
1+\varepsilon$, for any positive $\varepsilon\ll 1$.
First assume for simplicity that $D''\widetilde f=0$ without any error (i.e.\
that $\widetilde f$ is globally holomorphic). Then
$$
v_t=D''(\theta(\psi-t)\,\widetilde f)=\theta'(\psi-t)\,
d''\psi\wedge\widetilde f
\leqno(5.6)
$$
has support in the tubular domain $W_t=\{t<\psi<t+1\}$. At the same time, we adjust the functions $\eta=\eta_t$ and $\lambda=\lambda_t$ used in Prop.~3.11 to create enough convexity on $W_t$. For this, we take
$$
\eta_t=1-\delta\chi_t(\psi)\leqno(5.7)
$$
where $\chi_t:{}]-\infty,0]\to\bR$ is a negative smooth convex 
increasing function with the following properties:
{\plainitemindent=1.25cm
\smallskip
\item{(5.8~a)} $\chi_t(0)=0$~~ and~~ 
$\displaystyle\inf_{\tau\le 0}\chi_t(\tau)=-M_t>-\infty$,
\smallskip
\item{(5.8~b)} $0\le\chi'_t(\tau)\le{1\over 2}$~~ 
for $\tau\le 0$,
\smallskip
\item{(5.8~c)} $\chi'_t(\tau)=0$~~ for $\tau\in{}]-\infty,t-1]$,~~
$\,\chi'_t(\tau)>0$~~ for $\tau\in{}]t-1, 0]$,
\smallskip
\item{(5.8~d)} $\displaystyle\chi''_t(\tau)\ge {1-\varepsilon\over 4}$~~
for $\tau\in[t,t+1]$.
\smallskip
\item{(5.8~e)} $\displaystyle{\chi''_t(\tau)\over \chi'_t(\tau)^2}
\ge {2\delta\over\pi(1+\delta^2\tau^2)}$~~
for $\tau\in{}]t-1,0]$.
\smallskip
}
\noindent The function $\chi_t$ can be easily constructed by taking
$$
\chi''_t(\tau)={\delta\over 2\pi(1+\delta^2\tau^2)}\,\beta(\tau-t)+
{1-\varepsilon\over 4}\,\xi(\tau-t)\leqno(5.9)
$$
with support in $[t-1,0]$, where $\beta:\bR\to[0,1]$ is a smooth 
non decreasing function such that $\beta(\tau)=0$ for $\tau<-1$,
$\beta(\tau)=1$ for $\tau\ge 0$ and $\xi(\tau)=\beta(\tau)\beta(1-\tau)$.
We then have $\Supp(\beta)=[-1,+\infty[$, $\Supp(\xi)=[-1,2]$, and
$\xi(\tau)=\beta(\tau)$ on $[-1,0]$. On $[-1,0[$, we can take $\beta$ so 
small that $\int_{-1}^0\beta(\tau)\,d\tau<\varepsilon/2$. By symmetry we find
$\int_{-1}^2\xi(\tau)\,d\tau<1+\varepsilon$.
Clearly (5.8~a,c,d) hold if we adjust 
the integration constant so that $\chi_t'(t-1)=0$. Now,
the right hand side of (5.8) has a total integral on $]-\infty,0]$ 
that is less than $\smash{
{(1-\varepsilon)\over 4}(1+\varepsilon)+{1\over 2\pi}{\pi\over 2}<{1\over 2}}$,
thus (5.8~b) is satisfied. This implies that (5.8~e) holds at least on the
interval $[t,0]$. An integration by parts yields
$$
\chi'_t(\tau)\le\Big({\delta\over 2\pi}+{1-\varepsilon\over 4}\Big)
\widetilde\beta(\tau-t)\quad\hbox{on $[t-1,t]$},
$$
where $\widetilde\beta\ge 0$ is a primitive of $\beta$ vanishing at~$\tau=-1$. 
On $[t-1,t]$ we find
$$
{\chi''_t(\tau)\over \chi'_t(\tau)^2}\ge
{1-\varepsilon\over 4}\Big({\delta\over 2\pi}+{1-\varepsilon\over 4}\Big)^{-2}
{\beta(\tau-t)\over\widetilde\beta(\tau-t)^2}.\leqno(5.10)
$$
On $[-1,0]$ we have $\widetilde\beta(\tau)\le(1+\tau)\beta(\tau)\le\beta(\tau)$
and $\widetilde\beta(\tau)\le\int_{-1}^0\beta(\tau)\,d\tau<\varepsilon/2$, and
we see that the right hand side of (5.10) can be taken arbitrary large when 
$\varepsilon$ is small. Therefore (5.8~e) can also be achieved on $[t-1,t]$.

We now come back to the choice of $\eta_t$ and $\lambda_t$.
Since we have assumed $\psi\le 0$ at this step, 
we have by definition $\eta_t\ge 1$. Moreover, 
$d'\eta_t=-\delta\,\chi'_t(\psi)\,d'\psi$ and
$$
\ii\,d'd''\chi_t(\psi)=\ii\,\chi'_t(\psi)\,d'd''\psi+
\ii\,\chi''_t(\psi)\,d'\psi\wedge d''\psi,
\leqno(5.11)
$$
hence we see that
$$
\eqalign{
R_t&:=\eta_t\Big(\ii\,\Theta_{E,h}+\ii\,d'd''\psi\Big)-\ii\,d'd''\eta_t -
\lambda_t^{-1}\ii\,d'\eta_t\wedge d''\eta_t\cr
&\phantom{:}=\eta_t\Big(\ii\,\Theta_{E,h}+
(1+\delta\eta_t^{-1}\chi'_t(\psi))\,\ii\,d'd''\psi\Big)
+\Big(\delta\chi''_t(\psi) - \lambda_t^{-1}\delta^2\chi'_t(\psi)^2\Big)
\ii\,d'\psi\wedge d''\psi.\cr}
$$
The coefficient $(1+\delta\eta_t^{-1}\chi'_t(\psi))$ lies in 
$[1,1+{\delta\over 2}]$,
and our curvature assumption implies that the first term in the right hand
side is non negative. We take 
$$
\lambda_t=\pi(1+\delta^2\psi^2).
\leqno(5.12)
$$ 
By (5.8~e), this ensures that 
$\lambda_t^{-1}\delta^2\chi'_t(\psi)^2\le{1\over 2}\delta\chi''_t(\psi)$,
hence we get the crucial lower bounds
$$
\leqalignno{
&R_t\ge{1\over 2}\,\delta\chi''_t(\psi)\,\ii\,d'\psi\wedge d''\psi\ge 0
\quad\hbox{on $X_c$,}&(5.13)\cr
&R_t\ge{(1-\varepsilon)\delta\over 8}\,\ii\,d'\psi\wedge d''\psi
\quad\hbox{on $W_t=\{t<\psi<t+1\}$}.&(5.14)\cr}
$$
A standard calculation gives the formula
$\langle[\Lambda_\omega,\ii\,d'\psi\wedge d''\psi]v,v\rangle_{\omega,h}=
|(d''\psi)^*v|_{\omega,h}^2$ for any $(n,1)$-form $v$. Therefore, for every
$(n,0)$-form $u$ we have
$$
|\langle d''\psi\wedge u,v\rangle|^2=
|\langle u,(d''\psi)^*v\rangle|^2\le|u|^2|(d''\psi)^*v|^2=
|u|^2\langle[\Lambda_\omega,\ii\,d'\psi\wedge d''\psi]v,v\rangle.
$$
From this and (5.14) we infer that the curvature 
operator $B_t=[\Lambda_\omega,R_t]$
satisfies
$$
\langle B_t^{-1}(d''\psi\wedge u),d''\psi\wedge u\rangle=
|B_t^{-1/2}(d''\psi\wedge u)|^2\le
{8\over(1-\varepsilon)\delta}\,|u|^2\quad\hbox{on $W_t$}.
$$
In particular, since $|\theta'|\le 1+\varepsilon$, we see that the 
form $v_t=\theta'(\psi-t)\,d''\psi\wedge\widetilde f$ defined in (5.6) satisfies
$$
\langle B_t^{-1}v_t,v_t\rangle\le{8(1+\varepsilon)^2\over(1-\varepsilon)\delta}\,
|\widetilde f|^2.\leqno(5.15)
$$
Now, (5.8~b) implies $|\chi_t(\tau)|\le {1\over 2}|\tau|$, thus
$\eta_t=1-\delta\chi_t(\psi)\le 1-{1\over 2}\delta\psi$ and
$$\eta_t+\lambda_t\le 1-{\textstyle{1\over2}}\delta\psi+\pi(1+\delta^2\psi^2)
\le 4.21\,(1+\delta^2\psi^2)\leqno(5.16)
$$
(the optimal constant is ${\sqrt{5}+2\over 4}+\pi<4.21$). 
As $4.21\times 8\times(1+\varepsilon)^2/(1-\varepsilon)<34$ 
for $\varepsilon\ll 1$, Proposition 3.11 produces 
a solution $u_t$ such that $D''u_t=v_t$ on $X_c\ssm Y$ and
$$
\int_{X_c\ssm Y}(1+\delta^2\psi^2)^{-1}|u_t|^2_{\omega,h}e^{-\psi}\,
dV_{X,\omega}\le {34\over\delta}\int_{\{t<\psi<t+1\}}
|\widetilde f|^2_{\omega,h}e^{-\psi}dV_{X,\omega}.
$$
The function $F_t=\theta(\psi-t)\widetilde f -u_t$ is essentially the
extension we are looking for. For any $\alpha>0$ we have
$$
|F_t|^2\le (1+\alpha)\,|u_t|^2+(1+\alpha^{-1})\,
|\theta(\psi-t)|^2|\widetilde f|^2,
$$
hence
$$
\leqalignno{
\int_{X_c\ssm Y}(1&+\delta^2\psi^2)^{-1}(1+\alpha^2\psi^2)^{-(n-1)/2}
|F_t|^2_{\omega,h}\,e^{-\psi}\,dV_{X,\omega}\cr
&\le{34(1+\alpha)\over\delta}\int_{X_c\cap\{t<\psi<t+1\}}
|\widetilde f|^2_{\omega,h}e^{-\psi}dV_{X,\omega}&(5.17_1)\cr
&+(1+\alpha^{-1})\int_{X_c\cap\{\psi<t+1\}}(1+\delta^2\psi^2)^{-1}
(1+\alpha^2\psi^2)^{-(n-1)/2}\,
|\widetilde f|^2_{\omega,h}e^{-\psi}dV_{X,\omega}.&(5.17_2)\cr}
$$
The local integrability of $(1+|\psi|)^{-(n+1)}|\smash{\widetilde f}|^2e^{-\psi}$ 
asserted by Proposition~4.5~(c) and the Lebesgue dominated convergence theorem
imply that the  last integral $(5.17_2)$ converges to $0$ as~$t\to -\infty$. 
Therefore we obtain
$$
\limsup_{\alpha\to 0_+}\lim_{t\to-\infty}(5.17_1)+(5.17_2)
\le {34\over\delta}\int_{X_c\cap Y}|f|^2_{\omega,h}dV_{Y,\omega}[\psi].
$$
by definition of the measure in the right hand side (cf.\ (2.4)). 
Since $u_t$ is in $L^2$ with respect to the singular weight $e^{-\psi}$, it is
also locally $L^2$ with respect to a smooth weight. A standard lemma 
(cf.\ [Dem82, Lemme 6.9]) shows that the equation $D''u_t=v_t$ extends to~$X_c$,
and the hypoellipticity of $D''$ implies that $u_t$ is smooth on~$X_c$. As
$e^{-\psi}$ is non integrable along $Y$, we conclude that $u_t$ must vanish
on $X_c\cap Y$. Therefore $F_t$ is indeed a holomorphic extension of $f$ 
on $X_c$. By letting $t$ tend to $-\infty$ and then $c$ to $+\infty$, 
the uniformity of our $L^2$ inequalities implies that one can extract 
a weakly convergent sequence $F_{t_\nu}\to F$, such that 
$F\in H^0(X,K_X\otimes E)$ is an extension of $f$ and
$$
\int_X(1+\delta^2\psi^2)^{-1}|F|^2_{\omega,h}e^{-\psi}\,
dV_{X,\omega}\le {34\over\delta}\int_Y|f|^2_{\omega,h}dV_{Y,\omega}[\psi].
$$
The proof also works when $Y$ is singular, because the equations can be
considered on $X\ssm Y_{\rm sing}$, and $X_c\ssm Y_{\rm sing}$ is again 
complete K\"ahler for every~$c>0$. 

In case $\psi$ is no longer negative, we put
$\psi^+_A={1\over A}\log(1+e^{A\psi})$ and replace $\psi$ with
$$
\psi_A=\psi-\psi^+_A<0
$$
which converges to
$\psi-\psi_+=-\psi_-$ as $A\to+\infty$. We then solve $D''u_t=v_t$ with 
$v_t=D''(\theta(\psi_A-t)\widetilde f)$, and use the functions 
$\eta_t=1-\delta\chi_t(\psi_A)$ and $\lambda_t=\pi(1+\delta^2\psi_A^2)$
in the application of proposition~3.11. The expression of the curvature 
term $R_t$ becomes
$$
\eqalign{
R_{t,A}=\eta_t\Big(\ii\,\Theta_{E,h}&{}+\ii\,d'd''\psi+
\delta\eta_t^{-1}\chi'_t(\psi_A))\,\ii\,d'd''\psi_A\Big)\cr
&{}+\Big(\delta\chi''_t(\psi_A) - \lambda_t^{-1}\delta^2\chi'_t(\psi_A)^2\Big)
\ii\,d'\psi_A\wedge d''\psi_A.\cr}
$$
All bounds are then essentially the same, except that we have an additional
negative term $(...)\ii\,d'\psi\wedge d''\psi$ in $\ii\,d'd''\psi_A$, i.e.\
$$
\ii\,d'd''\psi_A={1\over 1+e^{A\psi}}\ii\,d'd''\psi - {Ae^{A\psi}\over
(1+e^{A\psi})^2}\ii\,d'\psi\wedge d''\psi.
$$
Because of this term, the first term $\eta_t(...)$ in $R_{t,A}$
is a priori no longer${}\ge 0$. However, this can be compensated by
adding an extra weight ${\delta\over 2}\psi^+_A$ to
the metric of $(E,h)$, so that the total weight is now
$\psi+{\delta\over 2}\psi^+_A$. The contribution of
the new weight to the term $\eta_t(...)$ in the modified $R_{t,A}$ reads
$$\eqalign{
\ii\,\Theta_{E,h}&+\Big(1+\delta\eta_t^{-1}\chi'_t(\psi_A){1\over 1+e^{A\psi}}+
{\delta\over 2}{e^{A\psi}\over 1+e^{A\psi}}\Big)\ii\,d'd''\psi\cr
&+\delta{Ae^{A\psi}\over(1+e^{A\psi})^2}
\Big({1\over 2}-\eta_t^{-1}\chi'_t(\psi_A)\Big)\ii\,d'\psi\wedge d''\psi,\cr}
$$
and as the coefficient of $\ii\,d'd''\psi$ still lies in $[1,1+\delta]$ (remember that $\eta_t\ge 1$ and $\chi'_t\le{1\over 2}$), we conclude that we have again
$$
\leqalignno{
&R_{t,A}\ge{1\over 2}\,\delta\chi''_t(\psi_A)\,\ii\,d'\psi_A\wedge d''\psi_A
\ge 0\quad\hbox{on $X_c$,}&(5.18)\cr
&R_{t,A}\ge{(1-\varepsilon)\delta\over 8}\,\ii\,d'\psi_A\wedge d''\psi_A
\quad\hbox{on $W_{t,A}=\{t<\psi_A<t+1\}$}.&(5.19)\cr}
$$
Since $\psi-\psi_A\to 0$ along $Y=\psi^{-1}(-\infty)$, it is easy to see that
the measures $dV_{Y,\omega}[\psi]$ and $dV_{Y,\omega}[\psi_A]$ coincide. After making those corrections, we get an extension $F_A$ such that
$$
\int_X(1+\delta^2\psi_A^2)^{-1}|F_A|^2_{\omega,h}e^{-\psi-{\delta\over 2}\psi^+_A}\,
dV_{X,\omega}\le {34\over\delta}\int_Y|f|^2_{\omega,h}dV_{Y,\omega}[\psi].
$$
By letting $A\to +\infty$ and extracting a limit $F_A\to F$, we get
$$
\int_X(1+\delta^2\psi_-^2)^{-1}|F|^2_{\omega,h}e^{-\psi-{\delta\over 2}\psi_+}\,
dV_{X,\omega}\le {34\over\delta}\int_Y|f|^2_{\omega,h}dV_{Y,\omega}[\psi],
$$
which is equivalent to the final bound given in Theorem~2.8.

The next point we have to justify is that unfortunately we cannot
expect $\smash D''\widetilde f\equiv 0$ as we assumed a priori (unless
we already know for some reason that a global holomorphic extension 
exists). In fact, we have to solve an equation 
$D''u=v_t:=D''(\theta(\psi_A-t)\widetilde f)$
with an extra term in the right hand side, namely
$$
D''u=v_t=v_t^{(1)}+v_t^{(2)},\quad
v_t^{(1)}=\theta'(\psi_A-t)d''\psi_A\wedge\widetilde f,\quad
v_t^{(2)}=\theta(\psi_A-t)D''\widetilde f.\leqno(5.20)
$$
The first term $v_t^{(1)}$ of $(5.20)$ has been already estimated, but we have 
to show that the second term $\smash{v_t^{(2)}}$ becomes ``negligible''
when we take limits as $t\to-\infty$. For this we solve (5.20) by means of
Prop.~3.12 instead of Prop.~3.11. We get an
approximate $L^2$ solution 
$D''u_{t,\varepsilon}=v_t-w_{t,\varepsilon}$, whence
$D''(\theta(\psi_A-t)\widetilde f - u_{t,\varepsilon})=w_{t,\varepsilon}$. Moreover,
this solution satisfies the $L^2$ estimate
$$
\leqalignno{
\Vert(\eta_t+\lambda_t)^{-1/2}u_{t,\varepsilon}\Vert^2&+
{1\over \varepsilon}\Vert w_{t,\varepsilon}\Vert^2\cr
&\le\int_{X_c\ssm Y}\langle (B_{t,A}+\varepsilon I)^{-1}v_t,v_t\rangle\,
e^{-\psi-{\delta\over 2}\psi^+_A}\,dV_{X,\omega}\cr
&\le(1+\alpha)\int_{X_c\ssm Y}\langle B_{t,A}^{-1}v_t^{(1)},v_t^{(1)}\rangle\,
e^{-\psi-{\delta\over 2}\psi^+_A}\,dV_{X,\omega}&(5.21_1)\cr
&+(1+\alpha^{-1})\varepsilon^{-1}
\int_{X_c\ssm Y}\langle v_t^{(2)},v_t^{(2)}\rangle\,
e^{-\psi-{\delta\over 2}\psi^+_A}\,dV_{X,\omega}&(5.21_2)\cr}
$$
for $\alpha>0$ arbitrary. The integral~$(5.21_1)$ is
bounded by means of $(5.17_i)$ -- or its analogue for the modified
weight $\psi+{\delta\over 2}\psi^+_A$ -- and the integral $(5.21_2)$ is 
equal to
$$
(1+\alpha^{-1})\varepsilon^{-1}
\int_{X_c\cap\{t<\psi<t+1\}}
|\theta(\psi_A-t)|^2|D''\widetilde f|^2\,e^{-\psi-{\delta\over 2}
\psi^+_A}\,dV_{X,\omega}.
\leqno(5.22)
$$
We now put all estimates together for the section
$F_{t,\varepsilon}=\theta(\psi_A-t)\widetilde f -u_{t,\varepsilon}$, which
satisfies
$$
|F_{t,\varepsilon}|^2\le (1+\alpha)\,|u_{t,\varepsilon}|^2+
(1+\alpha^{-1})\,|\theta(\psi_A-t)|^2|\widetilde f|^2.
$$
We get in this way from (5.16), $(5.17_i)$, $(5.21_i)$  and (5.22)
$$
\leqalignno{
\int_{X_c\ssm Y}(1&+\delta^2\psi_A^2)^{-1}(1+\alpha^2\psi_A^2)^{-(n-1)/2}
|F_{t,\varepsilon}|^2_{\omega,h}\,e^{-\psi-{\delta\over 2}\psi^+_A}\,dV_{X,\omega}\cr
&\kern150pt{}+{1\over\varepsilon}\int_{X_c\ssm Y}|w_{t,\varepsilon}|^2_{\omega,h}\,
e^{-\psi-{\delta\over 2}\psi^+_A}\,dV_{X,\omega}\cr
&\le{34(1+\alpha)^2\over\delta}
\int_{X_c\cap\{t<\psi_A<t+1\}}|\widetilde f|^2\,
e^{-\psi-{\delta\over 2}\psi^+_A}\,dV_{X,\omega}&(5.23_1)\cr
&+(1+\alpha^{-1})\int_{X_c\cap\{\psi_A<t+1\}}
(1+\delta^2\psi_A^2)^{-1}(1+\alpha^2\psi_A^2)^{-(n-1)/2}|\widetilde f|^2\,
e^{-\psi-{\delta\over 2}\psi^+_A}\,dV_{X,\omega}&(5.23_2)\cr
&+{(1+\alpha)(1+\alpha^{-1})\over\varepsilon}
\int_{X_c\cap\{\psi_A<t+1\}}|D''\widetilde f|^2\,
e^{-\psi-{\delta\over 2}\psi^+_A}\,dV_{X,\omega}.&(5.23_3)\cr}
$$
By Proposition~4.5~(c), the integral $(5.23_2)$ converges to $0$
as $t\to-\infty$. Since $|D''\widetilde f|^2\,e^{-\psi}$ is locally 
integrable on $X$, the last integral $(5.23_3)$ also converges to $0$ as 
$t\to-\infty$. We let $\varepsilon$ and $\alpha$ converge 
to $0$ afterwards, and extract a limit $F=
\lim_{\varepsilon\to 0}\lim_{t\to-\infty}F_{t,\varepsilon}$ as 
already explained, to recover the expected $L^2$ estimate.\smallskip

\noindent
{\bf (5.24) Final regularity argument.}
One remaining non trivial point is to check that we get smooth
solutions and that the resulting limit $F$ of 
$F_{t,\varepsilon}=\theta(\psi-t)\smash{\widetilde f}-u_{t,\varepsilon}$
is actually an extension of~$f$, one particular issue being that 
$F_{t,\varepsilon}$ 
is not exactly holomorphic (for the simplicity of notation we assume here
that $\psi\le 0$ since the difficulty is purely local near~$Y$, and thus
skip the $\psi_A$ approximation process in what follows).
For this, we apply
observation~5.1 and use a composition of blow-ups $\mu:\smash{\widehat X}
\to X$ such that the singularities of $\smash{\widehat\psi}=\psi\circ\mu$
are divisorial, given by some normal crossing divisor $\Delta=\sum
c_j\Delta_j$ in $\smash{\widehat X}$ whose support contains the
exceptional divisor of~$\psi$. If $g$ is a germ of section in
$\cO_{X,x}(K_X\otimes E)$, the section $\smash{\widehat g}=\mu^*g$
takes values in $\mu^*(K_X\otimes E)=K_{\widehat
X}\otimes(\cO_{\widehat X}(-\Delta')\otimes\mu^*E)$, where
$\Delta'=\sum c'_j\Delta_j$ has support in $|\Delta|$ and the $c'_j$
are non negative integers. As $\psi$ has log canonical singularities on~$X$,
we see by taking $g$ invertible that $\Delta-\Delta'$ has coefficients 
$c_j-c'_j\le 1$ on $\smash{\widehat X}$. Let us set
$$\Delta-\Delta'=\Delta^{(1)}+\Delta''$$
where $\Delta^{(1)}$ consists of the sum of components of multiplicity
$c_j-c'_j=1$, and $\Delta''$ is the sum of all other ones with
$c_j-c'_j<1$. The metric involved in our $L^2$ estimates
is $\mu^*h_\psi=\mu^*he^{-\smash{\widehat\psi}}$. When viewed as
a metric on $\smash{\widehat G}=\smash{K_{\widehat X}}\otimes\mu^*E
\otimes\cO(-\Delta'-\Delta^{(1)})$, it possesses a weight
$$
\widehat\psi_G:=\widehat\psi-\log|\sigma_{\Delta'}|^2-
\log|\sigma_{\Delta^{(1)}}|^2=\log|\sigma_{\Delta''}|^2\quad\mod~\cC^\infty,
$$
hence $\smash{\widehat\psi_G}$ is non singular at the generic point
of any component $\Delta_j$ of $\Delta^{(1)}$. Solving a 
$\smash{\dbar}$-equation in $\mu^*(K_X\otimes E)$ with respect
to the singular weight $e^{-\psi\circ\mu}=\smash{e^{-\widehat\psi}}$ amounts 
to solving the 
same $\dbar$-equation with values in $\smash{\widehat G}$
with respect to the weight $\smash{\widehat\psi_G}$. The standard 
ellipticity results imply that the
solutions $\widehat u_{t,\varepsilon}$ on $\smash{\widehat X}$ given by 
3.12 are smooth as sections of $\smash{\widehat G}$. This argument
shows that we can assume right away that $\psi$ has divisorial singularities,
and consider only the case $\mu=\Id_X$ and $\Delta=\Delta^{(1)}+\Delta''$, in
which case $Y=|\Delta^{(1)}|$ and $\cO_X(-\Delta^{(1)})=\cO_X(-Y)$. 
We thus simplify the notation by
removing all hats, and set $G=K_X\otimes E\otimes\cO_X(-Y)$. 
Notice that 
$$F_{t,\varepsilon}=\tau(\psi-t)\widetilde f-u_{t,\varepsilon}$$
is not exactly holomorphic. We have instead
$D''F_{t,\varepsilon}=w_{t,\varepsilon}$ where $w_{t,\varepsilon}$ 
is a section of $\Lambda^{0,1}T_X^*\otimes G$ that is smooth on the Zariski
open set $Y^\circ=Y\ssm\bigcup_{\Delta_k\not\subset Y}\Delta_k$, and its
$L^2$ norm with respect to $\psi_G=\psi-\log|\sigma_Y|^2$ satisfies
$\Vert w_{t,\varepsilon}\Vert=O(\varepsilon^{1/2})$. If $F_0$ is a fixed 
local extension
of $f$ near a generic point $x_0\in\Delta_j$, then $F_t-F_0$ is a local section of $G$ and \hbox{$D''(F_{t,\varepsilon}-F_0)=w_{t,\varepsilon}$}.
We can find, say by the standard H\"ormander $L^2$ estimates 
on a coordinate ball $B(x_0,r)$ ([H\"or65,66], see also [Kohn63,64]), 
a smooth $L^2$ section $s_{t,\varepsilon}$
of $G$ such that $\Vert s_{t,\varepsilon}\Vert=O(\varepsilon^{1/2})$
and $D''s_{t,\varepsilon}=w_{t,\varepsilon}$ on $B(x_0,r)$. Then
$F_{t,\varepsilon}-F_0-s_{t,\varepsilon}$ is a holomorphic
section of $G$ on~$B(x_0,r)$.
Its limit $F-F_0$ is a limit in $L^{2-\alpha}$ for every $\alpha>0$, thanks
to the H\"older inequality and the fact that $1+\delta^2\psi^2$ is in
$L^p$ for every $p>1$. Thus $F-F_0$ is a holomorphic section of $G$ on
$B(x_0,r)$, and so $F_{\restriction Y}=f$ on $Y^\circ$.\qed
\smallskip

\noindent{\bf Proof of Remark~2.9 (b).} By the technique of proof of
the regularization theorem [Dem92, Theorem~1.1], on any relatively compact
subset $X_c\compact X$, there are quasi-psh approximations 
$\varphi_\nu\downarrow\varphi$ of $\varphi$ with neat analytic singularities, 
such that the curvature estimate suffers only a small 
error${}\le 2^{-\nu}\omega$,
namely, for $h_\nu=e^{-\varphi_\nu}\le h$, we have
$$
\ii\,\Theta_{E,h_\nu}+\alpha\,\ii\,d'd''\psi_\nu\ge
-2^{-\nu}\omega\quad\hbox{on $X_c$},
$$
uniformly for $\alpha\in[1,1+\delta]$. As a consequence, the operator
$$B_{t,\nu}=\big[
\eta_t(\ii\,\Theta_{E,h_\nu}+\alpha\,\ii\,d'd''\psi_\nu
-\ii\,d'd''\eta_t-\lambda_t^{-1}\ii\,d'\eta_t\wedge d''\eta_t,
\Lambda_\omega]$$
has a slightly negative lower bound $-M_t2^{-\nu}$ by (5.8~a). This can
be absorbed by means of an additional positive term $\varepsilon I$ in
$B_{t,\nu}+\varepsilon I$, with $\varepsilon=O(2^{-\nu})$. Moreover the 
set of poles of $\varphi_\nu$ is an analytic set $P_\nu$ and we can work 
on the complete K\"ahler
manifold $X_c\ssm(Y\cup P_\nu)$ to avoid any singularities. Then
Prop.~3.12 provides an approximate solution $D''u_{t,\nu}\approx v_t$ 
with error $O(\varepsilon^{1/2})=O(2^{-\nu/2})$, satisfying the estimate
$$
\eqalign{
\int_{X_c\ssm Y}{\exp(-{\delta\over 2}{1\over A}\log(1+e^{A\psi}))\over
1+\delta^2\psi_A^2}\,
&|u_{t,\nu}|^2_{\omega,h_\nu}\,e^{-\psi}\,dV_{X,\omega}\le\cr
&{34\over\delta}\int_{X_c\cap\{t<\psi_A<t+1\}}
|\widetilde f|^2_{\omega,h_\nu}e^{-\psi_A}\,dV_{X,\omega}.\cr}
$$
The right hand side is uniformly bounded by a similar norm where
$h_\nu$ is replaced by $h=e^{-\varphi}$. The conclusion follows 
by letting $\nu$ converge to $+\infty$ (before doing anything else), 
and by extracting limits.\qed
\smallskip

\noindent{\bf Proof of Theorem~2.12.} The proof is essentially identical to
the proof of Theorem~2.8, we simply make a ``rescaling'': we replace 
$\psi$ by $m_p\psi$, $t$ by $m_pt$, $\delta$ by $\delta/m_p$, $\theta$ by 
$\theta(\tau)=\theta(\tau/m_p)$  and use Lemma~4.5 in its full generality.
Property (5.4) [namely the integrability of $|D''(\widetilde f)|^2e^{-m_p\psi}$]
still holds here since $f$ has local extensions $\widetilde f_j$ such that
the differences $\widetilde f_j-\widetilde f_k$ lie by construction 
in $\cI(m_p\psi)$.
In the final regularity argument 5.24, the vanishing of $f$
prescribed by $\cI(m_{p-1}\psi)$ subtracts a further divisor 
to $\Delta-\Delta'-\Delta_{m_{p-1}}$ where 
$\mu^*\cI(m_{p-1}\psi)=\cO(-\Delta_{m_{p-1}})$. The definition of jumps
leads to the fact that we only have to considers components of multiplicity
$1$ in that difference, and the rest of the argument is the same.\qed
\smallskip

\noindent{\bf Proof of Theorem 2.13.} The argument is identical to the proof
of Remark~2.9~(b), when we make the same rescaling, especially by
replacing $\psi$ with $m_p\psi$.\qed\smallskip

\noindent{\bf Proof of Theorem 2.14.} Since $X$ is compact, all coherent
cohomology groups involved are finite dimensional and Hausdorff for their
natural topology. We can assume $\psi<0$ by subtracting a constant
(and thus avoid the $\psi_A$ approximation process), and $X=X_c$.
\smallskip
\noindent
(a) We proceed by induction on~$p-\ell$. Let
$$f\in H^0(Y^{(m_p)},\cO_X(K_X\otimes E)\otimes\cI(m_\ell\psi)/\cI(m_p\psi)).$$
First assume $\ell=p-1$. Theorem 2.12 provides an extension $F$ of $f$ in
case $\cI(m_{p-1}\psi)$ is replaced by $\cI'(m_{p-1}\psi)$, but otherwise
the limiting $L^2$ integral of $f$ computed on $Y^{(m_p)}$ may diverge. 
The main idea, however, is that the integrals are still convergent when 
considered on each tube $\{t<\psi<t+1\}$, and this can be used to check
at least the qualitative part of the surjectivity theorem. In fact, we
apply the arguments used in the proof of Theorems~2.8--2.12, especially
estimates $(5.23_i)$, taken with $\alpha=1$, and $\psi=\psi_A$, $t$, 
$\theta(\tau)$, $\delta$ replaced respectively with $m_p\psi$, $m_pt$, 
$\theta(\tau/m_p)$, $\delta/m_p$; we also assume $\delta\le 1$ here.
This produces a $\cC^\infty$ extension
$F_{t,\varepsilon}=\theta(\psi-t)\,\widetilde f-u_{t,\varepsilon}$
such that $D''F_{t,\varepsilon}=w_{t,\varepsilon}$ and
$$
\leqalignno{
\int_{X}(1+m_p^2\psi^2&)^{-(n+1)/2}
|F_{t,\varepsilon}|^2_{\omega,h}\,e^{-m_p\psi}\,dV_{X,\omega}
+{1\over\varepsilon}\int_{X}|w_{t,\varepsilon}|^2_{\omega,h}\,
e^{-m_p\psi}\,dV_{X,\omega}\cr
\le{}&{136\over\delta}
\int_{\{t<\psi<t+1\}}|\widetilde f|^2\,
e^{-m_p\psi}\,dV_{X,\omega}&(5.25_1)\cr
&+2\int_{\{\psi<t+1\}}
(1+m_p\psi^2)^{-(n+1)/2}|\widetilde f|^2\,
e^{-m_p\psi}\,dV_{X,\omega}&(5.25_2)\cr
&+{4\over\varepsilon}
\int_{\{\psi<t+1\}}|D''\widetilde f|^2\,e^{-m_p\psi}\,dV_{X,\omega}.
&(5.25_3)\cr}
$$
This is true for all $t<0$ and $\varepsilon>0$, and the idea is to adjust
the choice of $\varepsilon$ as a function of~$t$. The integral
$(5.25_2)$ is uniformly bounded by Proposition~4.5~(c) and can be
disregarded. By construction the coefficients
$\smash{D''\widetilde f}$ can be expressed as a combination of functions
$\smash{\widetilde f_j-\widetilde f_j}$ in $\cI(m_p\varphi)$, thus by openness,
there exists $\alpha>0$ such that
$$
\int_X|D''\widetilde f|^2_{\omega,h}e^{-(m_p+\alpha)\psi}\,dV_{X,\omega}<+\infty.
$$
We therefore get an upper bound
$$
\int_{\{\psi<t+1\}}
|D''\widetilde f|^2_{\omega,h}e^{-m_p\psi}\,dV_{X,\omega}\le C\,e^{\alpha t}
\quad\hbox{as $t\to -\infty$}.\leqno(5.26)
$$
On the other hand, since 
$\int_X|\widetilde f|^2_{\omega,h}e^{-(m_p-\beta)\psi}\,dV_{X,\omega}<+\infty$
for every $\beta>0$, we conclude that
$$
\int_{\{t<\psi<t+1\}}
|\widetilde f|^2_{\omega,h}e^{-m_p\psi}\,dV_{X,\omega}\le C'(\beta)\,e^{-\beta t}
\quad\hbox{as $t\to -\infty$}.\leqno(5.27)
$$
A possible choice of $\varepsilon$ is to take $\varepsilon=e^{(\alpha+\beta)t}$
for some $\beta>0$. Then $(5.25_i)$ and (5.26--5.27) imply
$$
\int_X|w_{t,\varepsilon}|_{\omega,h}^2e^{-m_p\psi}dV_{X,\omega}
\le C''(\delta,\beta)e^{\alpha t}.
$$
We conclude that the error term $w_{t,\varepsilon}$ converges uniformly
to $0$ in $L^2$ norm as $t\to -\infty$ and $\varepsilon=e^{(\alpha+\beta)t}\to 0$.
The main term $F_{t,\varepsilon}$, however, is not under control, but we can
cope with this situation. By taking a principalization of 
the ideal defining the singularities of $\psi$, we can assume that these
singularities are divisorial, and thus that our solutions are smooth at 
the generic point of 
$Z_p=\Supp(\cI(m_{p-1}\psi)/\cI(m_p\psi))$. Observe that $w_{t,\varepsilon}$ is
a coboundary for the Dolbeault complex associated with the cohomology group
$H^0(X,\cO_X(K_X\otimes E)\otimes \cI(m_{p-1}\psi))$. Now, $X$ is compact
and \v{C}ech cohomology can be calculated on finite Stein coverings 
by spaces of \v{C}ech cocycles 
equipped with the topology given by $L^2$ norms with
respect to the weight $e^{-m_{p-1}\psi}$. We conclude via an isomorphism
between \v{C}ech cohomology and Dolbeault cohomology that there is a
smooth global section $s_{t,\varepsilon}$ of $\cC^\infty(K_X\otimes E)\otimes
\cI(m_{p-1}\psi)$ such that $D''s_{t,\varepsilon}=w_{t,\varepsilon}$ and
$$
\int_X|s_{t,\varepsilon}|^2_{\omega,h}e^{-m_{p-1}\psi}dV_{X,\omega}=O(e^{t\alpha}).
$$
On the other hand, on any coordinate ball $B=B(x_0,r)\subset X$, we can
apply the standard $L^2$ estimates of H\"ormander for bounded
pseudoconvex domains, and find another local solution 
$\widetilde s_{B,t,\varepsilon}$ such that
$D''\widetilde s_{B,t,\varepsilon}=w_{t,\varepsilon}$ and
$$
\int_{B(x_0,r)}|\widetilde s_{B,t,\varepsilon}|^2_{\omega,h}e^{-m_p\psi}dV_{X,\omega}=O(e^{t\alpha}).\leqno(5.28)
$$
The difference $s_{t,\varepsilon}-\widetilde s_{B,t,\varepsilon}$ is
a holomorphic section of $\cO_X(K_X\otimes E)\otimes\cI(m_{p-1}\psi)$
on $B(x_0,r)$ which converges to $0$ in $L^2$ norm, thus uniformly on
any smaller ball $B(x_0,r')$. Notice that $\widetilde s_{B,t,\varepsilon}$ has 
a vanishing order prescribed by $\cI(m_p\psi)$ by (5.28). Therefore
$F_{t,\varepsilon}-\widetilde s_{B,t,\varepsilon}$ is a local holomorphic section
of $\cO_X(K_X\otimes E)\otimes\cI(m_{p-1}\psi)$ that maps exactly to $f$.
By what we have seen
$$
F_{t,\varepsilon}-s_{t,\varepsilon}\in
H^0(X,\cO_X(K_X\otimes E)\otimes\cI(m_{p-1}\psi))
$$
is a global holomorphic section whose restriction converges uniformly
to $f$. Since $X$ is compact, we are dealing with finite dimensional
spaces of sections and thus the restriction morphism must be
surjective. The proof of the case $\ell=p-1$ is complete.

Now, assume that the result has been proved for $p-\ell=d$ and take 
$p-\ell=d+1$.
By reducing $f$ mod $\cI(m_{p-1}\psi)$, the induction hypothesis
provides an extension
$F'$ in $H^0(X,\cO_X(K_X\otimes E)\otimes\cI(m_\ell\psi))$ of
$$
f':=f~\mod~\cI(m_{p-1}\psi)\quad\hbox{in}\quad
H^0(Y^{(m_{p-1})},\cO_X(K_X\otimes E)\otimes\cI(m_\ell\psi)/\cI(m_{p-1}\psi)).
$$
Thus $f_p:=f-(F'~\mod~\cI(m_p\psi))$ defines a section
$$
f_p\in H^0(Y^{(m_p)},\cO_X(K_X\otimes E)\otimes\cI(m_{p-1}\psi)/\cI(m_p\psi)).
$$
The case $d=1$ provides an extension
$$
F_p\in H^0(X,\cO_X(K_X\otimes E)\otimes\cI(m_{p-1}\psi))
$$
and $F=F'+F_p$ is the extension we are looking for.
\medskip

\noindent
(b) Assume finally that $E$ is a line bundle and that $h$ is a
singular hermitian metric satisfying the curvature estimate in the
sense of currents only. We can reduce ourselves to the case when
$\psi$ has divisorial singularities. The regularization techniques of
[Dem15] (see Remark~3 before section~3, based on the solution of the
openness conjecture) produces a singular metric $h_\varepsilon$ withe
analytic singularities, such that the multiplier ideal sheaves
$\cI(he^{-m_\ell\psi})$ involved are unchanged when $h$ is replaced by
$h_\varepsilon$, and such that there is an arbitrary small loss in the
curvature.  We absorb this new adverse negative term by considering
again $B_t+\varepsilon I$, and by applying the same tricks as
above.\qed

\section{6}{References}
{\eightpoint\bigskip\parskip=2pt plus 1pt minus 1pt

\Bibitem [AN54]& Akizuki, Y., Nakano, S.& Note on Kodaira-Spencer's
proof of Lefschetz theorems& Proc.\ Jap.\ Acad.\ {\bf 30} (1954) 266--272&

\Bibitem [AV65]& Andreotti, A., Vesentini, E.& Carleman estimates for the
Laplace-Beltrami equation in complex manifolds& Publ.\ Math.\ I.H.E.S.\
{\bf 25} (1965) 81--130&

\Bibitem [Ber96]& Berndtsson, B.& The extension theorem of Ohsawa-Takegoshi
and the theorem of Donnelly-Fefferman& Ann.\ Inst.\ Fourier 
{\bf 14} (1996) 1087--1099&

\Bibitem [BL14]& Berndtsson, B., Lempert, L.& A proof of the 
Ohsawa-Takegoshi theorem with sharp estimates& arXiv:1407.4946&

\Bibitem [Blo13]& B{\l}ocki, Z.& Suita conjecture and the
Ohsawa-Takegoshi extension theorem& Invent.\ Math.\ {\bf 193} (2013),
no.~1, 149--158&

\Bibitem [Boc48]& Bochner, S.& Curvature and Betti numbers (I) and (II)&
Ann.\ of Math.\ {\bf 49} (1948), 379--390$\,$; {\bf 50} (1949), 77--93&

\Bibitem [Che11]& Chen, Bo-Yong& A simple proof of the Ohsawa-Takegoshi
extension theorem& arXiv: math.CV/ 1105.2430&

\Bibitem [Dem82]& Demailly, J.-P.& Estimations $L^2$ pour l'op\'erateur
$\dbar$ d'un fibr\'e vectoriel holomorphe semi-positif au dessus
d'une vari\'et\'e k\"ahl\'erienne compl\`ete& Ann.\ Sci.\ Ecole Norm.\
Sup.\ {\bf 15} (1982), 457--511&

\Bibitem [Dem92]& Demailly, J.-P.& Regularization of closed positive currents
and Intersection Theory& J.\ Alg.\ Geom.\ {\bf 1} (1992), 361--409&

\Bibitem [Dem97]& Demailly, J.-P.& On the Ohsawa-Takegoshi-Manivel
$L^2$ extension theorem& Proceedings of the Conference in honour of
the 85th birthday of Pierre Lelong, Paris, September 1997, \'ed.\
P.~Dolbeault, Progress in Mathematics, Birkh\"auser, Vol.~{\bf 188}
(2000), 47--82&

\Bibitem [Dem15]& Demailly, J.-P.& On the cohomology of pseudoeffective line
bundles& J.E. Forn{\ae}ss et al.\ (eds.), Complex Geometry and Dynamics, 
Abel Symposia 10, DOI 10.1007/978-3-319-20337-9\_4&

\Bibitem [Dem-X]& Demailly, J.-P.& Complex analytic and differential geometry&
self-published e-book, Institut Fourier, $455\,$pp, last version:
June 21, 2012&

\Bibitem [DHP13]& Demailly, J.-P., Hacon, Ch., P\u{a}un, M.&
Extension theorems, Non-vanishing and the existence of good minimal models&
Acta Math.\ {\bf 210} (2013), 203--259&

\Bibitem [DF83]& Donnelly H., Fefferman C.& $L^2$-cohomology and
index theorem for the Bergman metric& Ann.\ Math.\ {\bf 118} (1983),
593--618&

\Bibitem [DX84]& Donnelly H., Xavier F.& On the differential form spectrum 
of negatively curved Riemann manifolds& Amer.\ J.\ Math.\ {\bf 106}
(1984), 169--185&

\Bibitem [Gri66]& Griffiths, P.A.& The extension problem in complex analysis~II;
embeddings with positive normal bundle& Amer.\ J.~Math.\ {\bf 88} (1966),
366--446&

\Bibitem [Gri69]& Griffiths, P.A.& Hermitian differential geometry,
Chern classes and positive vector bundles& Global Analysis, papers in
honor of K.~Kodaira, Princeton Univ.\ Press, Princeton (1969),
181--251&

\Bibitem [GZ13]& Guan, Qi'an, Zhou, Xiangyu& Strong openness conjecture 
for plurisubharmonic functions& arXiv: math.CV/1311.3781&

\Bibitem [GZ15]& Guan, Qi'an, Zhou, Xiangyu& A solution of an $L^2$
extension problem with an optimal estimate and applications& Ann.\ of 
Math.\ (2) {\bf 181} (2015), no. 3, 1139--1208&

\Bibitem [Hir64]& Hironaka, H.& Resolution of singularities of an algebraic
variety over a field of characteristic zero& Ann.\ of Math.\ {\bf 79}
(1964), 109--326&

\Bibitem [H\"or65]& H\"ormander, L.& $L^2$ estimates and existence
theorems for the $\dbar$ operator& Acta Math.\ {\bf 113}
(1965), 89--152&

\Bibitem [H\"or66]& H\"ormander, L.& An introduction to Complex
Analysis in several variables& 1st edition, Elsevier Science Pub.,
New York, 1966, 3rd revised edition, North-Holland Math.\ library, Vol
7, Amsterdam (1990)&

\Bibitem [Kod53a]& Kodaira, K.& On cohomology groups of compact
analytic varieties with coefficients in some analytic faisceaux&
Proc.\ Nat.\ Acad.\ Sci.\ U.S.A.\ {\bf 39} (1953), 868--872&

\Bibitem [Kod53b]& Kodaira, K.& On a differential geometric method in the theory
of analytic stacks& Proc.\ Nat.\ Acad.\ Sci.\ U.S.A.\ {\bf 39} (1953),
1268--1273&

\Bibitem [Kohn63]& Kohn, J.J.& Harmonic integrals on strongly pseudo-convex
manifolds~I& Ann.\ Math.\ (2), {\bf 78} (1963), 206--213&

\Bibitem [Kohn64]& Kohn, J.J.& Harmonic integrals on strongly pseudo-convex
manifolds~II& Ann.\ Math.\ {\bf 79} (1964), 450--472&

\Bibitem [Kod54]& Kodaira, K.& On K\"ahler varieties of restricted type& Ann.\ 
of Math.\ {\bf 60} (1954), 28--48&

\Bibitem [Man93]& Manivel, L.& Un th\'eor\`eme de prolongement $L^2$ de
sections holomorphes d'un fibr\'e vectoriel& Math.\ Zeitschrift {\bf 212}
(1993), 107--122&

\Bibitem [MV07]& McNeal, J., Varolin, D.& Analytic inversion of adjunction: 
$L^2$ extension theorems with gain& Ann.\ Inst.\ Fourier (Grenoble)
{\bf 57} (2007), no.~3, 703--718&

\Bibitem[Nad89]& Nadel, A.M.& Multiplier ideal sheaves and
K\"ahler-Einstein metrics of positive scalar curvature& Proc.\ Nat.\
Acad.\ Sci.\ U.S.A. {\bf 86} (1989), 7299--7300~~~and~~
Annals of Math., {\bf 132} (1990), 549--596&

\Bibitem [Nak55]& Nakano, S.& On complex analytic vector bundles&
J.~Math.\ Soc.\ Japan {\bf 7} (1955), 1--12&

\Bibitem [Nak73]& Nakano, S.& Vanishing theorems for weakly $1$-complete
manifolds& Number Theory, Algebraic Geometry and Commutative
Algebra, in honor of Y.~Akizuki, Kinokuniya, Tokyo (1973), 169--179&

\Bibitem [Nak74]& Nakano, S.& Vanishing theorems for weakly $1$-complete
manifolds~II& Publ.\ R.I.M.S., Kyoto Univ.\ {\bf 10} (1974), 101--110&

\Bibitem [Ohs88]& Ohsawa, T.& On the extension of $L^2$ holomorphic
functions, II& Publ. RIMS, Kyoto Univ.\ {\bf 24} (1988), 265--275&

\Bibitem [Ohs94]& Ohsawa, T.& On the extension of $L^2$ holomorphic functions,
IV$\,$: A new density concept& Mabuchi, T.\ (ed.) et al.,
Geometry and analysis on complex manifolds. Festschrift for Professor
S.~Kobayashi's 60th birthday. Singapore: World Scientific,
(1994), 157--170&

\Bibitem [Ohs95]& Ohsawa, T.& On the extension of $L^2$ holomorphic
functions, III$\,$: negligible weights& Math.\ Zeitschrift {\bf 219}
(1995), 215--225&

\Bibitem [Ohs01]& Ohsawa, T.& On the extension of $L^2$ holomorphic
functions, V$\,$: Effects of generalization& Nagoya Math.\ J.\ {\bf 161}
(2001), 1--21, erratum: Nagoya Math.\ J.\ {\bf 163} (2001), 229&

\Bibitem [Ohs03]& Ohsawa, T.& On the extension of $L^2$ holomorphic
functions, VI$\,$: A limiting case& Explorations in complex and Riemannian
geometry, Contemp.\ Math., 332, Amer.\ Math.\ Soc., Providence, RI,
2003, 235--239&

\Bibitem [Ohs05]& Ohsawa, T.& $L^2$ extension theorems -- backgrounds and
a new result& Finite or infinite dimensional complex analysis and
applications, Kyushu University Press, Fukuoka, 2005, 261--274&

\Bibitem [OT87]& Ohsawa, T., Takegoshi, K.& On the extension of $L^2$
holomorphic functions& Math.\ Zeitschrift {\bf 195} (1987), 197--204&

\Bibitem [Pop05]& Popovici, D.& $L^2$ extension for jets of holomorphic sections of a Hermitian line bundle& Nagoya Math.\ J.\ {\bf 180} (2005), 1--34&

\Bibitem [Rie51]& Riemann, B.& Grundlagen f\"ur eine allgemeine Theorie 
der Functionen einer ver\"anderlichen complexen Gr\"osse&
Inauguraldissertation, G\"ottingen, 1851&

\Bibitem [Rie54]& Riemann, B.& Ueber die Hypothesen, welche der Geometrie 
zu Grunde liegen& Habilitationsschrift, 1854, Abhandlungen der K\"oniglichen
Gesellschaft der Wissenschaften zu Göttingen, 13 (1868)&

\Bibitem [Rie57]& Riemann, B.& Theorie der Abel'schen Functionen&
Journal f\"ur die reine und angewandte Mathematik, {\bf 54} (1857), 101--155&

\Bibitem [Var10]& Varolin, D.& Three variations on a theme in complex 
analytic geometry& Analytic and Algebraic Geometry, IAS/Park City
Math.\ Ser.\ {\bf 17}, Amer.\ Math.\ Soc., (2010), 183--294&

}
\bigskip\noindent
Jean-Pierre Demailly\\
Universit\'e Grenoble-Alpes\\
Institut Fourier, BP74\\
38400 Saint-Martin d'H\`eres, France\\
{\tt jean-pierre.demailly@ujf-grenoble.fr}

\bigskip\noindent
(version of October 14, 2015, printed on \today, \timeofday)
\end
\bye